\documentclass{amsart}
\usepackage{amsmath,amssymb,amsthm,amscd}
\newtheorem{formula}{}[section]
\newtheorem{proposition}[formula]{Proposition}
\newtheorem{corollary}[formula]{Corollary}
\newtheorem{lemma}[formula]{Lemma}
\newtheorem{theorem}[formula]{Theorem}
\theoremstyle{definition}
\newtheorem{definition}[formula]{Definition}
\newtheorem{example}[formula]{Example}
\theoremstyle{remark}
\newtheorem*{remark}{Remark}

\begin{document}
\title[Deformations of graded Lie algebras]
{Deformations of graded nilpotent Lie algebras and symplectic structures}
\author{Dmitri V. Millionschikov}
\thanks{Partially supported by
the Russian Foundation for Fundamental Research, grant no. 99-01-00090 and PAI-RUSSIE, dossier no. 04495UL}
\subjclass{17B30, 17B56, 17B70, 53D}
\address{Department of Mathematics and Mechanics, Moscow
State University, 119899 Moscow, RUSSIA}
\curraddr{Universit\'e Louis Pasteur, UFR de Math\'ematique et d'Informatique, 7 rue Ren\'e Descartes - 67084 Strasbourg Cedex (France)}
\email{million@mech.math.msu.su}

\begin{abstract}
We study symplectic structures on filiform Lie
algebras  -- nilpotent Lie algebras of the maximal length of the descending
central sequence. There are two basic examples of symplectic
${\mathbb Z}_{{>}0}$-graded filiform Lie algebras
defined by their basises
$e_1,\dots,e_{2k}$ and structure relations

1) $\mathfrak{m}_0(2k): \;[e_1,e_i]=e_{i{+}1}, \; i=2,\dots,2k{-}1.$

2) ${\mathcal V}_{2k}: \;[e_i,e_j]=(j{-}i)e_{i{+}j}, \; i{+}j\le 2k.$

Let $\mathfrak g$ be a symplectic filiform Lie algebra and $\dim \mathfrak{g}=2k \ge 12$. 
Then $\mathfrak g$ is isomorphic to some
${\mathbb Z}_{{>}0}$-filtered deformation either of $\mathfrak{m}_0(2k)$
or of ${\mathcal V}_{2k}$.
In the present article we classify ${\mathbb Z}_{{>}0}$-filtered deformations
of ${\mathcal V}_{n}$, i.e., Lie algebras with structure relations of the
following form:
$$[e_i,e_j]=(j{-}i)e_{i{+}j}+\sum_{l{=}1}c_{ij}^l e_{i{+}j{+}l}, \;\:i+j \le n$$
Namely we prove that  for $n \ge 16$ the moduli space ${\mathcal M}_n$ of these algebras
can be identified with the orbit space of the following
${\mathbb K}^*$-action on ${\mathbb K}^5$:
$$\alpha \star X=
(\alpha^{n{-}11}x_1,\alpha^{n{-}10}x_2,\dots,\alpha^{n{-}7}x_5), \alpha \in
{\mathbb K}^*, \; X \in {\mathbb K}^5.$$
For $n=2k$ the subspace
${\mathcal M}_{2k}^{sympl} \subset {\mathcal M}_{2k}$ of symplectic Lie algebras
is determined by equation $x_1=0$. A table with the structure constants
of symplecto-isomorphism classes in
${\mathcal M}_{2k}^{sympl}$ is presented.
\end{abstract}
\date{}

\maketitle

\section*{Introduction}
Nilmanifolds $M=G/\Gamma$ (compact homogeneous spaces of nilpotent
Lie groups $G$
over lattices $\Gamma$) are important examples of
symplectic manifolds that do not admit K\"ahler structures.
An interesting family of graded symplectic nilpotent Lie algebras ${\mathcal V}_{2k}$
(and corresponding family of nilmanifolds $M_n$)
was considered in \cite{BT1}, \cite{BT2}, \cite{Bu}. 
The finite dimensional Lie algebras ${\mathcal V}_n$, that are defined by the commutating relations 
$[e_i,e_j]=(j{-}i)e_{i{+}j}, \; i{+}j\le n$, "came" from the 
infinite dimensional Virasoro algebra and  they 
are examples
of so-called filiform Lie algebras --
nilpotent Lie algebras $\mathfrak{g}$ with the maximal
length $s{=}\dim \mathfrak{g}{-}1$ of the descending central sequence of
$\mathfrak{g}$. 
The study of filiform Lie algebras was started by
M.~Vergne in \cite{V1}, \cite{V2}. 

The classification of symplectic filiform Lie algebras
of dimensions $ \le 10$ was discussed in 
\cite{GJKh2}, \cite{GozeB}. The present paper is the continuation of \cite{Mill4},
where a criterion of the existence of a symplectic structure on a filiform Lie algebra
$\mathfrak g$ was proposed. In particular in a symplectic filiform 
$\mathfrak g=L^1 \mathfrak g$ one can find the ideal $L^2\mathfrak g$ of codimension $1$ 
such that the sequence of ideals $L^i \mathfrak g, i=1,\dots,2k$, where 
$L^i \mathfrak g=C^{i{-}1} \mathfrak g, i=3,\dots, n$  
($\{C^i \mathfrak g\}$ are the ideals of the descending central sequence $C$ of $\mathfrak g$)
determines a decreasing filtration $L$ of the Lie algebra $\mathfrak g$. 
The associated graded Lie algebra ${\rm gr}_L \mathfrak g$
is symplectic also. The graded filiform algebras of the type ${\rm gr}_L \mathfrak g$
were classified in \cite{Mill3}, \cite{Mill4}. There are two one-parameter families
of graded symplectic filiform Lie algebras of dimensions $8,10$. But
if ${\rm gr}_L \mathfrak g$ is a symplectic filiform Lie algebra with $\dim \mathfrak{g}=2k \ge 12$ 
then ${\rm gr}_L \mathfrak g$ is isomorphic either to $\mathfrak{m}_0(2k)$
or to ${\mathcal V}_{2k}$. In other words in dimensions $2k \ge 12$ one can obtain symplectic filiform
Lie algebras as  special deformations (that we call ${\mathbb Z}_{{>}0}$-filtered deformations) 
of two graded Lie algebras: $\mathfrak{m}_0(2k)$ and ${\mathcal V}_{2k}$.

We classify ${\mathbb Z}_{{>}0}$-filtered deformations
of ${\mathcal V}_{n}$, i.e., Lie algebras with the structure relations of the
following form:
$$[e_i,e_j]=(j{-}i)e_{i{+}j}+\sum_{l{=}1}c_{ij}^l e_{i{+}j{+}l}, \;\: i+j\le n.$$
We compute in the Section \ref{Computations} the space $H^2({\mathcal V}_{n}, {\mathcal V}_{n})$ for $n \ge 12$.
To the ${\mathbb Z}_{{>}0}$-filtered deformations corresponds the subspace
$\oplus_{{>}0}H^2_{(i)}({\mathcal V}_{n}, {\mathcal V}_{n})$.
The main theorem \ref{main} of the present article asserts that 
for $n \ge 16$ the moduli space ${\mathcal M}_n$ 
(i.e. the set of isomorphism classes) of these algebras
can be identified with the orbit space of the following
${\mathbb K}^*$-action on ${\mathbb K}^5=\oplus_{{>}0}H^2_{(i)}({\mathcal V}_{n}, {\mathcal V}_{n})$:
$$\alpha \star X=
(\alpha^{n{-}11}x_1,\alpha^{n{-}10}x_2,\dots,\alpha^{n{-}7}x_5), \alpha \in
{\mathbb K}^*, \; X \in {\mathbb K}^5,$$
where the coordinates $x_1,\dots, x_5$ of the space 
${\mathbb K}^5=\oplus_{{>}0}H^2_{(i)}({\mathcal V}_{n}, {\mathcal V}_{n})$ are defined by the
choice of the basic cocycles $\psi_{n,12{-}i}, i=1,\dots,5$. 

For $n=2k$ the subspace
${\mathcal M}_{2k}^{sympl} \subset {\mathcal M}_{2k}$ of symplectic Lie algebras
is determined by equation $x_1=0$. A table with structure constants
of symplecto-isomorphism classes in
${\mathcal M}_{2k}^{sympl}$ is presented in the Section 9.

\section{Invariant symplectic structures on Lie groups}

\begin{definition}  A Lie group $G$ is said to have a left-invariant 
symplectic structure if it has 
a left-invariant non-degenerate closed $2$-form $\omega$. 
\end{definition}

\begin{example}
$G$ is a two-dimensional abelian Lie group ${\mathbb R}^2$ with coordinates $x,y$ 
and $\omega=dx \wedge dy$.
\end{example}

\begin{example}
$G$ is a direct product ${\mathcal H}_3 \times {\mathbb R}$
of the Heisenberg group ${\mathcal H}_3$
of all matrices of the form
$$
   \left( \begin{array}{lcr}
   1 & x & z\\
   0 & 1 & y\\
   0 & 0 & 1\\
   \end{array} \right) , ~~~ x,y,z \in {\mathbb R},
$$
and abelian ${\mathbb R}$ (with coordinate $t$).
The invariant symplectic form $\omega$ is defined as
$$\omega = dx \wedge (dz - xdy)+dy \wedge dt.$$
\end{example}

\begin{theorem}[B.-Y.~Chu, ~\cite{Chu}]
Any semisimle Lie group has no left-invariant symplectic structure.
\end{theorem}

\begin{theorem}[B.-Y.~Chu, ~\cite{Chu}]
A connected unimodular Lie group admitting a left-invariant
symplectic structure must to be solvable.
\end{theorem}

One can obtain examples of left-invariant symplectic structures  
in the framework of Kirillov's orbit method.

Let us consider the coadjoint action $Ad^*$ of a Lie group $G$ on
the dual $\mathfrak{g}^{*}$ of its Lie algebra $\mathfrak{g}$:
$$(Ad^*(g)f)(X)=f(Ad(g^{{-}1})X), \; \forall X \in {\mathfrak g},$$
where $g \in G, f \in \mathfrak{g}^*$. The orbit $O_f$ of this action
has a homogeneous symplectic structure $\omega_{O_f}$ such that
$$\pi^*_f(\omega_{O_f})=df,$$ 
where $\pi_f$ denotes the natural mapping
$\pi_f: G \to O_f, \pi_f(g)=Ad^*(g)f$. Let the stabilizer
$G_f=\{g \in G| Ad^*(g)f=f\}$ of $f \in \mathfrak{g}^*$ be
a normal subgroup of $G$ then
the orbit $O_f$ can be identified with the quotient group $G/G_f$ and
the corresponding symplectic structure $\omega_{O_f}$ is left
$G/G_f$-invariant.
If $G$ is a nilpotent Lie group then nilpotent ones are $G_f$ and $G/G_f$.

\begin{definition}
A skew-symmetric non-degenerate bilinear form $\omega$ on the Lie algebra $\mathfrak{g}$
is called symplectic if it closed, i.e.
$$\omega([x,y],z)+\omega([y,z],x)+\omega([z,y],x)=0, \; \forall x,y,z \in \mathfrak{g}.$$
\end{definition}

If $\omega_G$ is a left-invariant symplectic form on $G$,
then $\omega_G$ defines a symplectic structure $\omega_{\mathfrak g}$ 
on the Lie algebra $\mathfrak g$ of $G$, and conversely any
symplectic form $\omega_{\mathfrak g}$ of ${\mathfrak g}$ defines 
a left-invariant symplectic structure on $G$.

\begin{lemma}[A.~M\'edina, P.~Revoy ~\cite{MR}]
Let ${\mathfrak g}$ be a symplectic Lie algebra with
non-trivial center $Z({\mathfrak g})$ and $I$ be a one-dimensional
subspace in $Z({\mathfrak g})$ and
$I^{\omega}$ its  symplectic complement with respect
to $\omega$.
Then one can consider two
following exact sequences of Lie algebras and their homomorphisms
\begin{equation}
\begin{split}
0 \to I \to I^{\omega}  \to I^{\omega}/I \to 0\\
0 \to I^{\omega}  \to {\mathfrak g}
\to {\mathfrak g}/I^{\omega} {\cong} I \to 0
\end{split}
\end{equation}
where $I^{\omega}/I$ is symplectic Lie algebra
(the restriction of $\omega$ to $I^{\omega}$ defines 
symplectic form $\tilde \omega$ on the quotient-algebra $I^{\omega}/I$).
\end{lemma}
In other words $2k$-dimensional symplectic Lie algebra ${\mathfrak g}$
with non-trivial center 
can be obtained from $2k{-}2$-dimensional symplectic
$I^{\omega}/I$ by means of two consequtive operations:

1) one-dimensional central extension of $I^{\omega}/I$ by $I$;

2) semidirect product of $I^{\omega}$ and one-dimensional
${\mathfrak g}/I^{\omega} \cong I$.

The combination of these two operations
is called the double extension of symplectic Lie algebra $I^{\omega}/I$.

\begin{theorem}[A.~M\'edina, P.~Revoy ~\cite{MR}]
Symplectic nilpotent Lie algebras can be obtained
by means of the sequence of consecutive double extensions starting
with trivial Lie algebra of zero dimension.
\end{theorem}

The even-dimensional nilpotent Lie algebras are classified for
dimensions $2k \le 6$ (~\cite{Mor}).
The classification of symplectic $6$-dimensional nilpotent Lie algebras 
based on this classification
was done in ~\cite{GozeB}.

\section{Filiform Lie algebras}

The sequence of ideals of a Lie algebra $\mathfrak{g}$
$$C^1\mathfrak{g}=\mathfrak{g} \; \supset \;
C^2\mathfrak{g}=[\mathfrak{g},\mathfrak{g}] \; \supset \; \dots
\; \supset \;
C^{k}\mathfrak{g}=[\mathfrak{g},C^{k-1}\mathfrak{g}] \; \supset
\; \dots$$
is called the descending central sequence of $\mathfrak{g}$.

A Lie algebra $\mathfrak{g}$ is called nilpotent if
there exists $s$ such that:
$$C^{s+1}\mathfrak{g}=[\mathfrak{g}, C^{s}\mathfrak{g}]=0,
\quad C^{s}\mathfrak{g} \: \ne 0.$$
The natural number $s$ is called the nil-index of the nilpotent Lie algebra $\mathfrak{g}$, or 
$\mathfrak{g}$ is  called  $s$-step nilpotent Lie algebra.

Let $\mathfrak{g}$ be a Lie algebra. We call a set $F$ of subspaces
$$\mathfrak{g} \supset \dots \supset F^{i}\mathfrak{g} \supset
F^{i{+}1}\mathfrak{g} \supset \dots \qquad (i \in \mathbb Z)
$$
a decreasing filtration $F$ of $\mathfrak{g}$ if $F$ is compatible with the Lie structure
$$[F^{k}\mathfrak{g},F^{l}\mathfrak{g}] \subset F^{k+l}\mathfrak{g},\; \forall k,l \in \mathbb Z.$$

Let $\mathfrak{g}$ be a filtered Lie algebra.  
A graded Lie algebra 
$${\rm gr}_F\mathfrak{g}=\bigoplus_{k=1} ({\rm gr}_F\mathfrak{g})_k, \;\;
({\rm gr}_F\mathfrak{g})_k=F^{k}\mathfrak{g}/F^{k{+}1}\mathfrak{g}$$
is called the associated graded Lie algebra
${\rm gr}_F\mathfrak{g}$.

The ideals $C^{k}\mathfrak{g}$ of the descending central sequence define
a decreasing filtration $C$ of the Lie algebra
$\mathfrak{g}$ 
$$C^{1}\mathfrak{g}=\mathfrak{g} \supset C^{2}\mathfrak{g} \supset
\dots \supset C^{k}\mathfrak{g} \supset \dots; \qquad
[C^{k}\mathfrak{g},C^{l}\mathfrak{g}] \subset C^{k+l}\mathfrak{g}.$$
One can consider the associated graded Lie algebra
${\rm gr}_C\mathfrak{g}$.

The finite filtration $C$ of a nilpotent Lie algebra $\mathfrak g$ is called 
the canonical filtration of
a nilpotent Lie algebra $\mathfrak g$. 

\begin{proposition}
Let $\mathfrak{g}$ be a $n$-dimensional nilpotent Lie algebra.
Then for its nil-index we have the estimate $s \le n-1$.
\end{proposition}

\begin{definition}
A nilpotent $n$-dimensional Lie algebra $\mathfrak{g}$ is  called
filiform Lie algebra if it has the nil-index $s=n-1$. 
\end{definition}

\begin{example}
The Lie algebra $\mathfrak{m}_0(n)$ is defined 
by its basis $e_1, e_2, \dots, e_n$
with commutating relations:
$$ [e_1,e_i]=e_{i+1}, \; \forall \; 2\le i \le n{-}1.$$
\end{example}

\begin{remark}
We will omit in the sequel trivial
commutating relations $[e_i,e_j]=0$ in the definitions of Lie algebras.
\end{remark}

\begin{example}
The Lie algebra $\mathfrak{m}_2(n)$ is
defined by its basis $e_1, e_2, \dots, e_n$
and commutating relations:
$$
[e_1, e_i ]=e_{i+1}, \quad 2 \le i \le n{-}1; \quad \quad
[e_2, e_j ]=e_{j+2}, \quad  3 \le j \le n{-}2. 
$$
\end{example}

\begin{example}
Let us define the algebra $L_k$ as the infinite-dimensional Lie algebra of
polynomial vector fields
on the real line ${\mathbb R}^1$ with a zero in $x=0$ of
order not less then $k+1$.

The algebra   $L_k$ can be defined by its infinite basis and  commutating relations
$$e_i=x^{i+1}\frac{d}{dx}, \; i \in {\mathbb N},\; i \ge k, \quad \quad
[e_i,e_j]= (j-i)e_{i{+}j}, \; \forall \;i,j \in {\mathbb N}.$$

One can consider the $n$-dimensional quotient algebra ${\mathcal V}_n=L_1/L_{n{+}1}$. 
\end{example}

The Lie algebras $\mathfrak{m}_0(n)$, $\mathfrak{m}_2(n)$, ${\mathcal V}_n$
considered above are filiform Lie algebras.

\begin{proposition}
Let $\mathfrak{g}$ be a filiform Lie algebra and
${\rm gr}_C\mathfrak{g}= \oplus_i ({\rm gr}_C\mathfrak{g})_i$ is
the corresponding associated (with respect to the  canonical filtration $C$)
graded Lie algebra. Then 
$$
\dim ({\rm gr}_C\mathfrak{g})_1=2, \quad
\dim ({\rm gr}_C\mathfrak{g})_2=\dots=
\dim ({\rm gr}_C\mathfrak{g})_{n{-}1}=1.
$$
\end{proposition}

We have the following isomorphisms of graded Lie algebras:
$$ {\rm gr}_C\mathfrak{m}_2(n) \cong
{\rm gr}_C{\mathcal V}_n \cong {\rm gr}_C\mathfrak{m}_0(n) \cong
\mathfrak{m}_0(n).$$

\begin{theorem}[M. Vergne \cite{V2}]
\label{V_ne1}
Let $\mathfrak{g}=\oplus_{\alpha}\mathfrak{g}_{\alpha}$ be a
graded $n$-dimensional filiform Lie algebra and 
\begin{equation}
\label{can_grad}
\dim \mathfrak{g}_1=2, \quad
\dim \mathfrak{g}_2=\dots=
\dim \mathfrak{g}_{n{-}1}=1.
\end{equation}
then  

1) if $n=2k+1$, then $\mathfrak{g}$ is isomorphic to
$\mathfrak{m}_0(2k+1)$;

2) if $n=2k$, then $\mathfrak{g}$ is isomorphic either to
$\mathfrak{m}_0(2k)$ or to the Lie algebra $\mathfrak{m}_1(2k)$,
defined by its basis $e_1, \dots, e_{2k}$ 
and commutating relations:
$$
[e_1, e_i ]=e_{i+1}, \; i=2, \dots, 2k{-}1; \quad \quad
[e_j, e_{2k{+}1{-}j} ]=(-1)^{j{+}1}e_{2k}, \quad j=2, \dots, k.
$$
\end{theorem}

\begin{remark}
In the settings of the Theorem \ref{V_ne1}
the gradings of the algebras $\mathfrak{m}_0(n)$, $\mathfrak{m}_1(n)$
are defined as
$\mathfrak{g}_1=Span( e_1, e_2 ) $,
$\mathfrak{g}_i=Span( e_{i{+}1} ), \: i=2, \dots, n{-}1.$
\end{remark}

\begin{corollary}[M. Vergne \cite{V2}]
Let $\mathfrak g$ be a filiform Lie algebra. Then one can choose a so-called
adapted basis $e_1, e_2, \dots, e_n$ in $\mathfrak g$:
\begin{equation}
[e_1, e_i]=e_{i+1}, i{=}2, {\dots}, n{-}1; \quad
[e_i, e_j]= \left\{\begin{array}{r}
   \sum \limits_{k{=}0}^{n{-}i{-}j}c_{ij}^{i{+}j{+}k} e_{i{+}j{+}k}, \; i+j \le n;\\
   (-1)^i \alpha e_n, \; i+j=n+1 ; \\
   0, \hspace{2.76em} \; i+j > n+1; \\
   2 \le i < j \le n.
   \end{array} \right. 
\end{equation}
where $\alpha =0$ if $n$ is odd number.
\end{corollary}

\section{Symplectic filiform Lie algebras: filtrations and gradings} 

\begin{definition}
A Lie algebra $\mathfrak{g}$ is called symplectic if it admits at least one symplectic structure.
\end{definition}

\begin{lemma}[~\cite{Mill4}]
Let $\mathfrak g$ be an $2k$-dimensional symplectic filiform Lie algebra,
then
\begin{equation}
\label{nec_cond}
{\rm gr}_C\mathfrak{g} \cong \mathfrak{m}_0(2k).
\end{equation}
\end{lemma}

In other words: 
let $\mathfrak g$ be a symplectic filiform Lie algebra and       
$e_1, e_2, \dots, e_{2k}$ be some adapted basis, then $[e_i, e_j]=0, \; i+j=2k+1$. 

\begin{proposition}
A fixed  adapted basis of $\mathfrak g$ $e_1, e_2, \dots, e_{2k}$  such that 
$[e_i, e_j]=0, \; i+j=2k+1$ defines a non-canonical 
filtration $L$ of $\mathfrak g$ :
\begin{equation}
\begin{split}
\mathfrak g=L^1 \mathfrak g \supset L^2 \mathfrak g \supset \dots \supset
L^{2k} \mathfrak g \supset \{0\},\\
L^j\mathfrak{g}{=}Span ( e_{j}, \dots, e_{2k} ),
\, j=1,{\dots},2k.\\ 
\end{split}
\end{equation}
\end{proposition}

\begin{remark} For homogeneous components of associated graded ${\rm gr}_L\mathfrak{g}$ we have  
$$
\dim ({\rm gr}_L\mathfrak{g})_1=
\dim ({\rm gr}_L\mathfrak{g})_2=\dots=
\dim ({\rm gr}_L\mathfrak{g})_{2k}=1.
$$
\end{remark}

\begin{proposition}[~\cite{Mill4}]
Let $\mathfrak{g}$ be a symplectic filiform Lie algebra with a fixed
adapted basis.
Then the corresponding associated graded Lie algebra 
${\rm gr}_L\mathfrak{g}= \oplus_i ({\rm gr}_L\mathfrak{g})_i$
is symplectic also.
\end{proposition}

\begin{remark}
It was shown in ~\cite{Mill4} that the previuos condition is necessary
but not sufficient condition.
\end{remark}

\begin{theorem}[D. Millionschikov, ~\cite{Mill4}]
Let $\mathfrak{g}=\bigoplus_{\alpha=1}^{2k}\mathfrak{g}_{\alpha}$ be
a real graded symplectic Lie algebra such that
\begin{equation}
\label{gradfil}
\dim \mathfrak{g}_{\alpha} = 1, \; 1 \le \alpha \le 2k; \quad \quad
[\mathfrak{g}_1,\mathfrak{g}_{\alpha}]=\mathfrak{g}_{{\alpha}{+}1}, \;
2 \le \alpha \le 2k-1.
\end{equation}
Then $\mathfrak{g}$ is isomorphic to the one and
only one Lie algebra  from the following list :

$$
\begin{tabular}{|c|c|c|c|}
\hline
&&&\\[-10pt]
dim & algebra & commutating relations & symplectic form\\
&&&\\[-10pt]
\hline
&&&\\[-10pt]
$4$ & $\mathfrak{m}_0(4)$ &
$[e_1, e_2]=e_{3}, [e_1, e_3]=e_{4}$ & $e^1{\wedge}e^4-e^2{\wedge}e^3$\\
&&&\\[-10pt]
\hline
&&&\\[-10pt]
$6$ & $\mathfrak{m}_0(6)$ &
$[e_1, e_i]=e_{i{+}1}, i=2,\dots,5$ & $e^1{\wedge}e^6-e^2{\wedge}e^5+e^3{\wedge}e^4$\\
&&&\\[-10pt]
\cline{2-4}
&&&\\[-10pt]
 & ${\mathcal V}_6$ &
$[e_i, e_j]=(j-i)e_{i{+}j},  i {+} j \le 6$ & $5e^1{\wedge}e^6+3e^2{\wedge}e^5+e^3{\wedge}e^4$\\
&&&\\[-10pt]
\hline
&&&\\[-10pt]
$8$ & $\mathfrak{m}_0(8)$ &
$[e_1, e_i]=e_{i{+}1}, i=2,\dots,7$ & $e^1{\wedge}e^8-e^2{\wedge}e^7+e^3{\wedge}e^6-e^4{\wedge}e^5$\\
&&&\\[-10pt]
\cline{2-4}
&&&\\[-10pt]
& 
${\mathfrak g}_{8,\alpha}$ 
 & 
 \begin{tabular}{c}
 (1)\\
 $\alpha \ne {-}\frac{5}{2}, {-}2, {-}\frac{1}{2},\frac{1}{2}$\\ \end{tabular} &
\begin{tabular}{c}
$e^1{\wedge}e^8{+}\frac{2\alpha^2{+}3\alpha{-}2}{2\alpha{+}5}e^2{\wedge}e^7{+}$\\
${+}\frac{2\alpha{+}2}{2\alpha{+}5}e^3{\wedge}e^6+\frac{3}{2\alpha{+}5}e^4{\wedge}e^5$\\
\end{tabular}\\
&&&\\[-10pt]
\hline
&&&\\[-10pt]
$10$ & $\mathfrak{m}_0(10)$ &
$[e_1, e_i]=e_{i{+}1}, i=2,\dots,9$ & 
$e^1{\wedge}e^{10}{-}e^2{\wedge}e^9{+}e^3{\wedge}e^8{-}e^4{\wedge}e^7{+}e^5{\wedge}e^6$\\
&&&\\[-10pt]
\cline{2-4}
&&&\\[-10pt]
& 
${\mathfrak g}_{10,\alpha}$ 
 & 
 \begin{tabular}{c}
 (1),(2)\\
 $\alpha \ne {-}\frac{5}{2}, {-}\frac{1}{4}, {-}1, {-}3, \alpha_1, \alpha_2;$\\ 
 $\alpha_1,\alpha_2 \in {\mathbb R},\: 2\alpha_1^3{+}2\alpha_1^2{+}3=0,$\\
$4\alpha_2^3{+}8\alpha_2^2{-}8\alpha_2{-}21=0$\\ 
 \end{tabular} &
\begin{tabular}{c}
$e^1{\wedge}e^{10}{+}\frac{2\alpha^3{+}2\alpha^2{+}3}{2(\alpha^2{+}4\alpha{+}3)}e^2 {\wedge} e^9{+}$\\
$+\frac{4\alpha^3{+}8\alpha^2{-}8\alpha{-}21}
{2(\alpha^2{+}4\alpha{+}3)(2\alpha{+}5)}e^3 {\wedge} e^8+$\\
$+\frac{3(2\alpha^2{+}4\alpha{+}5)}
{2(\alpha^2{+}4\alpha{+}3)(2\alpha{+}5)}e^4 {\wedge} e^7+$\\
$+\frac{3(4\alpha{+}1)}{2(\alpha^2{+}4\alpha{+}3)(2\alpha{+}5)}e^5 {\wedge} e^6$\\
\end{tabular}\\
&&&\\[-10pt]
\hline
&&&\\[-10pt]
$2k\ge 12$ & $\mathfrak{m}_0(2k)$ &
$[e_1, e_i]=e_{i{+}1}, i=2,\dots,2k{-}1$ & $\frac{1}{2}\sum_{i{+}j{=}2k{+}1} ({-}1)^{i{+}1}e^i{\wedge}e^j$\\
&&&\\[-10pt]
\cline{2-4}
&&&\\[-10pt]
 & ${\mathcal V}_{2k}$ &
$[e_i, e_j]=(j-i)e_{i{+}j},  i {+} j \le 2k$ & $\frac{1}{2}\sum_{i{+}j{=}2k{+}1}(j{-}i)e^i{\wedge}e^j$
\\[2pt]
\hline
\end{tabular}
$$

$$
\begin{tabular}{|c|c|}
\hline
&\\[-10pt]
& commutating relations for $\mathfrak{g}_{8, \alpha}$ and $\mathfrak{g}_{10, \alpha}$\\
&\\[-10pt]
\hline
&\\[-10pt]
(1) & \begin{tabular}{c}
$[e_1,e_2]=e_3, \; [e_1,e_3]=e_4, \;
[e_1,e_4]=e_5, \; [e_1,e_5]= e_6, \; [e_1,e_6]=e_7,$ \\
$[e_2,e_3]=(2{+}\alpha) e_5, \; [e_2,e_4]=(2{+}\alpha) e_6, \;
[e_2,e_5]=(1{+}\alpha)e_7, \; [e_3,e_4]= e_7$\\
$[e_1,e_7]=e_8, \quad [e_2,e_6]=\alpha e_8, \quad [e_3,e_5]= e_8$ \\
\end{tabular}\\
&\\[-10pt]
\hline
&\\[-10pt]
(2) & \begin{tabular}{c}
$[e_1,e_8]=e_9, \; [e_2,e_7]=\frac{2\alpha^2{+}3\alpha{-}2}{2\alpha{+}5}e_9, \;
[e_3,e_6]= \frac{2\alpha{+}2}{2\alpha{+}5}e_9, \;
[e_4,e_5]= \frac{3}{2\alpha{+}5}e_9$\\
$[e_1,e_9]=e_{10}, \;
[e_2,e_8]=\frac{2\alpha^2{+}\alpha{-}1}{2\alpha{+}5}e_{10}, \;
[e_3,e_7]= \frac{2\alpha{-}1}{2\alpha{+}5}e_{10}, \;
[e_4,e_6]= \frac{3}{2\alpha{+}5}e_{10}$\\
\end{tabular}\\[4pt]
\hline
\end{tabular}$$

\begin{remark}
In the fourth column of the table
we give only one variant of possible symplectic structure.  
\end{remark}

\end{theorem}

\begin{corollary}
Let $\mathfrak g$ be a symplectic filiform Lie algebra of dimension $2k \ge 12$ then 
one can choose a basis $e_1, \dots, e_{2k}$ in $\mathfrak g$ such that  
the corresponding commutating relations will be either 
\begin{equation}
\begin{split}
[e_1, e_i]=e_{i+1}+\sum \limits_{l{=}1}^{2k{-}i{-}1}c_{1i}^{l} e_{i{+}l{+}1}, \; i = 2, \dots, 2k{-}1; \\
[e_i, e_j]= \sum \limits_{l{=}1}^{2k{-}i{-}j}c_{ij}^{l} e_{i{+}j{+}l}, \; i+j \le 2k;
   2 \le i < j \le 2k;   
\end{split}
\end{equation}
or 
\begin{equation}
\begin{split}
[e_i, e_j]= (j-i)e_{i{+}j}+\sum \limits_{l{=}1}^{2k{-}i{-}j}c_{ij}^l e_{i{+}j{+}l}, \;
 i+j \le 2k;
   1 \le i < j \le 2k;   
\end{split}
\end{equation}

\end{corollary}

\begin{remark} 
The one-parameter family $\mathfrak{g}_{8, \alpha}$ was considered 
in ~\cite{GozeB} as well as corresponding symplectic form $\omega_{9}(\alpha)$. In \cite{GJKh2}
symplectic (over $\mathbb C$) low-dimensional ($\dim \mathfrak g \le 10$) 
filiform Lie algebras were classified (but this article contains  some mistakes).
\end{remark}

\section{Lie algebra cohomology}
\label{cohomology}
Let $\mathfrak{g}$ be a Lie algebra over ${\mathbb K}$ and
$\rho: \mathfrak{g} \to \mathfrak{gl}(V)$ its linear representation
(or in other words $V$ is a $\mathfrak{g}$-module).
We denote by $C^q(\mathfrak{g},V)$
the space of $q$-linear skew-symmetric mappings of $\mathfrak{g}$ into
$V$. Then one can consider an algebraic complex:

$$
\begin{CD}
V @>{d_0}>>
C^1(\mathfrak{g}, V) @>{d_1}>> C^2(\mathfrak{g}, V) @>{d_2}>>
\dots @>{d_{q{-}1}}>> C^q(\mathfrak{g}, V) @>{d_q}>> \dots
\end{CD}
$$
where the differential $d_q$ is defined by:

\begin{equation}
\begin{split}
(d_q f)(X_1, \dots, X_{q{+}1})=
\sum_{i{=}1}^{q{+}1}(-1)^{i{+}1}
\rho(X_i)(f(X_1, \dots, \hat X_i, \dots, X_{q{+}1}))+\\
+ \sum_{1{\le}i{<}j{\le}q{+}1}(-1)^{i{+}j{-}1}
f([X_i,X_j],X_1, \dots, \hat X_i, \dots, \hat X_j, \dots, X_{q{+}1}).
\end{split}
\end{equation}

The cohomology of the complex $(C^*(\mathfrak{g}, V), d)$ is called
the cohomology of the Lie algebra $\mathfrak{g}$
with coefficients in the representation $\rho: \mathfrak{g} \to V$.

In this article we will consider two main examples:

1) $V= {\mathbb K}$ and $\rho: \mathfrak{g} \to {\mathbb K}$ is trivial;

2) $V= \mathfrak{g}$ and $\rho=ad: \mathfrak{g} \to \mathfrak{g}$ is
the adjoint representation of $\mathfrak{g}$.

The cohomology of $(C^*(\mathfrak{g}, {\mathbb K}), d)$ (the first example)
is called the
cohomology with trivial coefficients of the Lie algebra
$\mathfrak{g}$ and is denoted by $H^*(\mathfrak{g})$.
Also we fixe
the notation $H^*(\mathfrak{g},\mathfrak{g})$ for the cohomology
of $\mathfrak{g}$ with coefficients in the adjoint representation.

One can remark that
$d_1: C^1(\mathfrak{g}, {\mathbb K}) \rightarrow
C^2(\mathfrak{g}, {\mathbb K})$ of the $(C^*(\mathfrak{g}, {\mathbb K}), d)$
is the dual mapping to the Lie bracket
$[ \, , ]: \Lambda^2 \mathfrak{g} \to \mathfrak{g}$. 
Moreover the condition $d^2=0$ is equivalent to the Jacobi identity for $[,]$.

Let $\mathfrak{g}=\oplus_{\alpha}\mathfrak{g}_{\alpha}$ be a
${\mathbb Z}$-graded Lie algebra
and $V=\oplus_{\beta} V_{\beta}$ is a ${\mathbb Z}$-graded
$\mathfrak{g}$-module, i.e.,
$\mathfrak{g}_{\alpha}V_{\beta} \subset V_{\alpha {+} \beta}.$
Then the complex $(C^*(\mathfrak{g}, V), d)$
can be equipped with the ${\mathbb Z}$-grading
$C^q(\mathfrak{g},V) =
\bigoplus_{\mu} C^q_{(\mu)}(\mathfrak{g},V)$, where
a $V$-valued $q$-form $c$ belongs to
$C^q_{(\mu)}(\mathfrak{g},V)$  
iff for
$X_1 \in \mathfrak{g}_{\alpha_1}, 
\dots,  X_q \in \mathfrak{g}_{\alpha_q}$ we have
$$c(X_1,\dots,X_q) \in
V_{\alpha_1{+}\alpha_2{+}\dots{+}\alpha_q{+}\mu}.$$ 

This grading is compatible with the differential  $d$ and
hence we have ${\mathbb Z}$-grading in cohomology:
$$
H^{q} (\mathfrak{g},V)= \bigoplus_{\mu \in {\mathbb Z}}
H^{q}_{(\mu)} (\mathfrak{g},V).
$$

\begin{remark}
The trivial $\mathfrak{g}$-module ${\mathbb K}$ has only one non-trivial
homogeneous component ${\mathbb K}={\mathbb K}_0$.
\end{remark}

\begin{example}
Let $\mathfrak{g}$ be a Lie algebra with the basis
$e_1, e_2, \dots, e_n$ and commutating relations
$$[e_i,e_j]= c_{ij}e_{i{+}j},  i+j \le n.$$

Let us consider the dual basis $e^1, e^2, \dots, e^n$.
One can introduce a grading (that we will call the weight)
of $\Lambda^*(\mathfrak{g}^*)=C^*(\mathfrak{g},{\mathbb K})$:
$$\Lambda^* (\mathfrak{g}^*)=
\bigoplus_{\lambda{=}1}^{n(n{+}1)/2}
\Lambda^*_{(\lambda)} (\mathfrak{g}^*),$$
where a subspace 
$\Lambda^{q}_{(\lambda)} (\mathfrak{g}^*)$ is spanned by $q$-forms
$\{ e^{i_1} {\wedge} \dots {\wedge} e^{i_q}, \;
i_1{+}\dots{+}i_q {=} \lambda \}$. 
For instance a monomial 
$e^{i_1} \wedge \dots \wedge e^{i_q}$ has the degree $q$ and the weight 
$\lambda=i_1{+}\dots{+}i_q$.

The complex $(C^*(\mathfrak{g}, \mathfrak{g}), d)$ is
${\mathbb Z}$-graded: 
$$C^*(\mathfrak{g}, \mathfrak{g})=
\bigoplus_{\mu \in {\mathbb Z}}C^*_{(\mu)}(\mathfrak{g}, \mathfrak{g}),$$
where $C^q_{(\mu)}(\mathfrak{g}, \mathfrak{g})$ is spanned
by monomials
$\{ e_l \otimes e^{i_1} {\wedge} \dots {\wedge} e^{i_q}, \;
i_1{+}\dots{+}i_q{+}\mu =l \}$. 
\end{example}

Now we consider a filtered Lie algebra $\mathfrak{g}$ with 
a finite decreasing filtration $F$.
One can define a decreasing filtration $\tilde F$ of
$\Lambda^* ({\mathfrak{g}}^*)$.
$$ \tilde F^{\mu} \Lambda^p (\mathfrak{g}^*)= \left\{ \omega
\in \Lambda^p (\mathfrak{g}^*)
 \; | \; \omega(F^{\alpha_1} \mathfrak{g} \wedge \dots \wedge F^{\alpha_p}\mathfrak{g})=0, \:  
 \alpha_1{+} \dots {+} \alpha_p {+} \mu \ge 0 \right\}.$$

\begin{example}
Let $\mathfrak{g}$ be a Lie algebra with the basis
$e_1, e_2, \dots, e_n$ and commutating relations
$$
[e_i,e_j]= \sum\limits_{k=0}^{n{-}i{-}j} c^{k}_{ij}e_{i{+}j{+}k},  i+j \le n.
$$
As it was remarked above the corresponding filtration $L$ ($L^k=Span(e_k,\dots,e_n)$) of $\mathfrak g$ can be defined.
The associated graded Lie algebra ${\rm gr}_L \mathfrak g$ has the following
structure relations
$$
[e_i,e_j]= c_{ij}^{0} e_{i{+}j},  i+j \le n.
$$

Let us consider the dual basis $e^1, e^2, \dots, e^n$. Then 
$\tilde L^{\mu}\Lambda^{p} (\mathfrak{g}^*)$ is spanned by $p$-monomials of weights
less or equal to $-\mu$, i.e. by
$p$-forms $e^{i_1} {\wedge} \dots {\wedge} e^{i_p}$ such that $i_1{+}\dots{+}i_p {\le} -\mu$. For instance
$$\tilde L^{{-}5}\Lambda^{2} (\mathfrak{g}^*)=Span(e^1 \wedge e^2, e^1 \wedge e^3, e^1 \wedge e^4, e^2 \wedge e^3).$$  
\end{example}

\begin{remark}
One can consider the spectral sequence $E_r$ that corresponds to the filtration $\tilde F$ of
the complex $\Lambda^* ({\mathfrak{g}}^*)$. We have an isomorphism (see \cite{V2} for example)
$$E^{p,q}_1 = H^{p{+}q}_{(-p)}({\rm gr}_F \mathfrak g).$$  
\end{remark}

\begin{theorem}[~\cite{Mill4}]
\label{criterion}
Let $\mathfrak g$ be a filiform Lie algebra such that 
${\rm gr}_C \mathfrak{g} \cong \mathfrak{m}_0(2k)$ and ${\rm gr}_L \mathfrak{g}$ is symplectic.

Then the Lie algebra $\mathfrak g$ is symplectic if and only if some homogeneous symplectic
class $[\omega_{2k{+}1}] \in E_1^{{-}2k{-}1, 2k{+}3}=H^{2}_{(2k{+}1)}({\rm gr}_L \mathfrak g)$ survives 
to the term $E_{\infty}$. 
\end{theorem}

\section{${\mathbb Z}_{{>}0}$-filtered deformations and $H^2(\mathfrak{g},\mathfrak{g})$}
\label{Nijenhuis-Richardson}

In this section we recall some definitions from the Nijenhuis-Richardson
deformation theory (see ~\cite{NR}).

\begin{definition}
Let $\mathfrak g$ be a Lie algebra with a Lie bracket $[,]$
and $\Psi: {\mathfrak g} \otimes {\mathfrak g} \to {\mathfrak g}$
is a skew-symmetric bilinear map. $\Psi$ is called a deformation
of $[,]$ iff $[,]'=[,]+\Psi$ is a Lie bracket on the vector space ${\mathfrak g}$.
\end{definition}

The Jacobi identity for bracket $[,]'=[,]+\Psi$ 
$$\left[[x,y]',z\right]'+\left[[y,z]',x\right]'+\left[[z,x]',y\right]'=0$$
is equivalent to the so-called deformation equation:
\begin{equation}
\begin{split}
\Psi([x,y],z){+}\Psi([y,z],x){+}\Psi([z,x],y){+}
\left[\Psi(x,y),z \right]{+}\left[\Psi(y,z),x \right]{+}
\left[\Psi(y,z),x \right]{+}\\
\\ {+}\Psi(\Psi(x,y),z){+}\Psi(\Psi(y,z),x){+}
\Psi(\Psi(z,x),y)=0.
\end{split}
\end{equation}
The first six terms can be rewrited in the form
$d\Psi(x,y,z)$ where
$d:C^2(\mathfrak g,\mathfrak g) \to C^3(\mathfrak g,\mathfrak g)$ is
the differential of the complex $(C^*(\mathfrak g,\mathfrak g),d)$.

Finally we have
\begin{equation}
\label{mdef}
d\Psi+\frac{1}{2}[\Psi,\Psi]=0,
\end{equation}
where $[,]$ denotes a symmetric bilinear function
$[,]: C^2(\mathfrak g,\mathfrak g)\times
C^2(\mathfrak g,\mathfrak g) \to C^3(\mathfrak g,\mathfrak g)$
\begin{equation}
\begin{split}
[\Psi,\tilde \Psi](x,y,z)=\Psi(\tilde \Psi(x,y),z)+\Psi(\tilde \Psi(y,z),x)+
\Psi(\tilde \Psi(z,x),y)+ \\
+\tilde \Psi(\Psi(x,y),z)+\tilde \Psi(\Psi(y,z),x)+\tilde \Psi(\Psi(z,x),y).
\end{split}
\end{equation}

The last definition can be generalised in terms of 
Nijenhuis-Richardson bracket in $C^*(\mathfrak g,\mathfrak g)$:
$$
[,]: C^p(\mathfrak g,\mathfrak g)\times
C^q(\mathfrak g,\mathfrak g) \to
C^{p{+}q{-}1}(\mathfrak g,\mathfrak g).$$
Namely, for $\alpha \in C^p(\mathfrak g,\mathfrak g)$ and
$\beta \in C^q(\mathfrak g,\mathfrak g)$ one can define
$[\alpha,\beta] \in C^{p{+}q{-}1}(\mathfrak g,\mathfrak g)$:
\begin{equation}
\begin{split}
[\alpha{,}\beta](\xi_1,{\small \dots},\xi_{p{+}q{-}1}){=}
\sum\limits_{1{\le}i_1{<}{\small \dots}{<}i_q{\le}p{+}q{-}1}
\alpha(\beta(\xi_{i_1}{,}{\small \dots}{,}\xi_{i_q})\xi_1{,}{\small \dots}{,}
\hat \xi_{i_1}{,}{\small \dots}{,}
\hat \xi_{i_q}{,} {\small \dots}{,}\xi_{p{+}q{-}1}){+}\\
+({-}1)^{pq{+}p{+}q}\sum\limits_{1{\le}j_1{<}{\dots}{<}j_p{\le}p{+}q{-}1}
\beta(\alpha(\xi_{j_1}{,}{\small \dots}{,}\xi_{j_p})\xi_1{,}{\small \dots}{,}
\hat \xi_{j_1}{,}{\small \dots}{,}
\hat \xi_{j_q}{,} {\small \dots}{,}\xi_{p{+}q{-}1}).
\end{split}
\end{equation}

The Nijenhuis-Richardson bracket defines a Lie superalgebra structure
in $C^*(\mathfrak g,\mathfrak g)$, i.e., if
$\alpha \in C^p(\mathfrak g,\mathfrak g)$,
$\beta \in C^q(\mathfrak g,\mathfrak g)$ and
$\gamma \in C^r(\mathfrak g,\mathfrak g)$ then

\begin{equation}
\label{superalg}
\begin{split}
1)\quad [\alpha,\beta]={-}({-}1)^{(p{-}1)(q{-}1)}[\beta,\alpha];
\hspace{30mm}\\
2) \; ({-}1)^{(p{-}1)(q{-}1)} \left[[\alpha,\beta],\gamma\right]{+}
({-}1)^{(q{-}1)(r{-}1)}\left[[\beta,\gamma],\alpha\right]{+}
({-}1)^{(r{-}1)(p{-}1)}\left[[\gamma,\alpha],\beta\right]= 0.
\end{split}
\end{equation}
Also we have the following important property:
$$d[\alpha,\beta]=[d\alpha,\beta]+({-}1)^{p}[\alpha,d\beta].$$
Thus the Nijenhuis-Richardson bracket defines a Lie superalgebra structure
in cohomology $H^*(\mathfrak g,\mathfrak g)$, i.e., the set of
bilinear functions
$$
[,]: H^p(\mathfrak g,\mathfrak g)\times
H^q(\mathfrak g,\mathfrak g) \to
H^{p{+}q{-}1}(\mathfrak g,\mathfrak g)$$
with the properties (\ref{superalg}).

\begin{proposition}
Let $\mathfrak{g}=
\oplus_{\alpha} \mathfrak{g}_{\alpha}$
be a ${\mathbb Z}$-graded Lie algebra. 
Then the ${\mathbb Z}$-grading of $C^*(\mathfrak{g},\mathfrak{g})$
and $H^*(\mathfrak{g},\mathfrak{g})$
is compatible with the Nijenhuis-Richardson bracket:
\begin{equation}
\begin{split}
[,]: C^p_{(\mu)}(\mathfrak{g},\mathfrak{g}) \times
C^q_{(\nu)}(\mathfrak{g},\mathfrak{g}) \longrightarrow
C^{p{+}q{-}1}_{(\mu+\nu)}(\mathfrak{g},\mathfrak{g})\\
[,]: H^p_{(\mu)}(\mathfrak{g},\mathfrak{g}) \times
H^q_{(\nu)}(\mathfrak{g},\mathfrak{g}) \longrightarrow
H^{p{+}q{-}1}_{(\mu+\nu)}(\mathfrak{g},\mathfrak{g})
\end{split}
\end{equation}
\end{proposition}

\begin{definition}
A deformation $\Psi$ of a ${\mathbb Z}_{{>}0}$-graded Lie algebra
$\mathfrak{g}=\oplus_{\alpha > 0} \mathfrak{g}_{\alpha}$ is called ${\mathbb Z}_{{>}0}$-filtered 
if the following condition
holds:
$$\forall X_1 \in \mathfrak{g}_{\alpha_1},
\dots, \forall X_q \in \mathfrak{g}_{\alpha_q}, 
\quad \Psi(X_1,\dots,X_q) \in
\bigoplus_{\alpha {>}
\alpha_1{+}\alpha_2{+}\dots{+}\alpha_q} V_{\alpha}.$$
Or in other words:
\begin{equation}
\label{kdef}
\Psi=\Psi_{1}+\Psi_{2}+\dots+\Psi_i+\dots, \quad \Psi_i \in C^2_{(i)}(\mathfrak{g},\mathfrak{g}), 
\:i=1,2,\dots 
\end{equation}
\end{definition}

Decomposing $\Psi=\Psi_{1}+\Psi_2+\dots$ in the deformation equation (\ref{mdef}) and
comparing terms with the same grading we came to the following system
of equations on homogeneous components $\Psi_l$: 
\begin{equation}
\label{sysdef}
\begin{split}
d\Psi_1=0, \:
d\Psi_{2}+\frac{1}{2}[\Psi_{1},\Psi_{1}]=0, \quad
d\Psi_{3}+[\Psi_{1},\Psi_{2}]=0, \dots, \\
d\Psi_{i}+\frac{1}{2}\sum_{m{+}l{=}i}[\Psi_{m},\Psi_{l}]=0, \dots
\hspace{19mm}
\end{split}
\end{equation}
This is the well-known system of equations for one-parametric deformation,
in our case one can associate to $\Psi$ the following one-parametric
deformation $\Psi_t$:
$$
\Psi_t=t\Psi_{1}+t^{2}\Psi_{2}+\dots+t^i\Psi_i+\dots
$$

1) the first equality of (\ref{sysdef}) shows that
$\Psi_1$ is a cocycle, we denote by $\bar \Psi_1$ its cohomology class
in $H^2_{(1)}(\mathfrak{g},\mathfrak{g})$.

2) the second one shows that the Nijenhuis-Richarson product
$[\bar \Psi_1,\bar \Psi_1]$ defines a trivial element in $H^3_{(1)}(\mathfrak{g},\mathfrak{g})$. 
If $[\bar \Psi_1,\bar \Psi_1]=0$ then $\Psi_2$ is determined
not uniquely but up to some closed element in $Z^2_{(2)}(\mathfrak{g},\mathfrak{g})$. 

3) hence one can find $\Psi_3$
iff $[\Psi_1,\Psi_2]$ is trivial in
$H^3_{(3)}(\mathfrak{g},\mathfrak{g})$ for some choice of $\Psi_2$.

The subset in $H^3_{(3)}(\mathfrak{g},\mathfrak{g})$ formed by 
elements ${-}[\Psi_1,\Psi_2]$, where $\Psi_2$ is a solution in
$C^2_{(2)}(\mathfrak{g},\mathfrak{g})$
of the equation $d\Psi_{2}+\frac{1}{2}[\Psi_{1},\Psi_{1}]=0$ is
called a triple Massey product $<\bar \Psi_1,\bar \Psi_1,\bar \Psi_1>$ 
and it is defined iff $[\bar \Psi_1,\bar \Psi_1]=0$.
\begin{definition}
A Massey product of the $k$-th order
$<\bar \Psi_1,\bar \Psi_1,\dots,\bar \Psi_1>$
is a subset in $H^3_{(k)}(\mathfrak{g},\mathfrak{g})$ formed by
classes $-\frac{1}{2}\sum_{m{+}l{=}k}[\Psi_{m},\Psi_{l}]$
where $\Psi_{j}, j{=}2,{\dots},k{-}1$ are solutions of the first $k{-}1$
equations of the system (\ref{sysdef}) and $\Psi_1$ is a representative
of $\bar \Psi_1$. A Massey product $<\bar \Psi_1,\bar \Psi_1,\dots,\bar \Psi_1>$
is called trivial if it contains zero.
\end{definition}
\begin{remark}
A Massey product $<\bar \Psi_1,\bar \Psi_1,\dots,\bar \Psi_1>$ of
the the $k$-th order
is defined iff all Massey products
$<\bar \Psi_1,\bar \Psi_1,\dots,\bar \Psi_1>$
of orders less than $k$ are trivial. One can easily show
that it does not depend on the choice of $\Psi_1$ in $\bar \Psi_1$.
A Lie product $[\bar \Psi_1,\bar \Psi_1]$
is called a Massey product of the second order.
\end{remark}
\begin{proposition}
Let be $\Psi_1$ an element of $H^2_{(1)}(\mathfrak{g},\mathfrak{g})$.
One can construct a deformation $\Psi$ of $\mathfrak{g}$ with a first term
equal to $\Psi_1 \in H^2_{(1)}(\mathfrak{g},\mathfrak{g})$ if and only
if all Massey products
$<\bar \Psi_1,\bar \Psi_1,\dots,\bar \Psi_1>$ are trivial.
\end{proposition}

Let $\Psi$ and $\tilde\Psi$ be two deformations of a Lie algebra
$\mathfrak{g}$.
The question is: whether they define non-isomorphic Lie algebras
or not ? Or does there exist a non-degenerate linear transformation
$\varphi:\mathfrak{g} \to \mathfrak{g}$ such that
\begin{equation}
\varphi([x,y]{+}\tilde \Psi(x,y))=
\left[\varphi(x),\varphi(y)\right]+
\Psi(\varphi(x),\varphi(y))
\end{equation}

Let us consider an non-degenerate linear operator
of the form $\varphi=id+\varphi_1, \: \varphi_{1} \in C^1_{(1)}(\mathfrak{g}, \mathfrak{g})$.
We have for the first term  
$$
\varphi^{{-}1}\left(
\left[\varphi(x),\varphi(y)\right]{+}
\Psi(\varphi(x),\varphi(y))\right)=[x,y]+(\Psi_1+d\varphi_1)(x,y)+\dots
$$

\begin{proposition} 
\label{prpdef}
Let $\Psi$ be some ${\mathbb Z}_{{>}0}$-filtered deformation of $[,]$ and $\Phi_1$ an arbitrary cocycle
in $Z^2_{(1)}(\mathfrak{g}, \mathfrak{g})$ cohomologous to $\Psi_1$.
Then the Lie algebra with $[,]+\Psi$ is isomorphic to some filtered deformation
$[,]+\tilde \Psi, \tilde \Psi_1=\Phi_1$.
\end{proposition}

\section{Computations}
\label{Computations}
In this section we compute $H^2({\mathcal V}_n, {\mathcal V}_n)$ for $n \ge 12$.

Let $\mathfrak{g}=
\oplus_{\alpha >0} \mathfrak{g}_{\alpha}$
be a ${\mathbb Z}_{{>}0}$-graded Lie algebra. One can define
a decreasing filtration ${\mathcal F}$ of
$(C^*(\mathfrak{g},\mathfrak{g}),d)$:
$$
{\mathcal F}^0 C^*(\mathfrak{g},\mathfrak{g})\supset \dots
\supset {\mathcal F}^q C^*(\mathfrak{g},\mathfrak{g})
\supset {\mathcal F}^{q{+}1} C^*(\mathfrak{g},\mathfrak{g}) \supset \dots
$$
where the subspace ${\mathcal F}^q C^{p{+}q}(\mathfrak{g},\mathfrak{g})$
is spanned by $p{+}q$-forms $c$ in $C^{p{+}q}(\mathfrak{g},\mathfrak{g})$
such that   
$$
c(X_1,\dots,X_{p{+}q}) \in \bigoplus_{\alpha=q} {\mathfrak g}_{\alpha} \:
\forall X_1,\dots,X_{p{+}q} \in \mathfrak{g}.
$$

The filtration ${\mathcal F}$ is compatible with $d$.

Let us consider the corresponding spectral sequence $E_r^{p,q}$:
\begin{proposition}
$E_1^{p,q}={\mathfrak g}_q \otimes H^{p{+}q}({\mathfrak g})$.
\end{proposition}

We have the following natural isomorphisms: 
\begin{equation}
\begin{split}
C^{p{+}q}({\mathfrak g},{\mathfrak g})
= {\mathfrak g} \otimes \Lambda^{p{+}q}({\mathfrak g}^*)\\
E_0^{p,q}={\mathcal F}^q C^{p{+}q}(\mathfrak{g},\mathfrak{g})/
{\mathcal F}^{q{+}1} C^{p{+}q}(\mathfrak{g},\mathfrak{g})
={\mathfrak g}_q \otimes \Lambda^{p{+}q}({\mathfrak g}^*).
\end{split}
\end{equation}

Now the proof follows from the formula for the $d_0^{p,q}: E_0^{p,q} \to E_0^{p{+}1,q}$:
$$d_0(X \otimes f)=X\otimes df,$$
where $X \in \mathfrak{g},
f \in \Lambda^{p{+}q}({\mathfrak g}^*)$ and $df$ is the standart
differential of the cochain complex of $\mathfrak{g}$ with trivial
coefficients.

\begin{theorem}
Let $n \ge 12$ then

1) $\dim H^0({\mathcal V}_n, {\mathcal V}_n)=
\dim H^0_{(n)}({\mathcal V}_n, {\mathcal V}_n)=1$;

2) $\dim H^1({\mathcal V}_n, {\mathcal V}_n)=4$,
namely
$$\dim H^1_{(\mu)}({\mathcal V}_n, {\mathcal V}_n)
= \left\{\begin{array}{r}
   1, \quad \mu=0,n{-}4, n{-}3, n{-}2;
   \hspace{0mm}\\
   0, \hspace{6mm} otherwise. \hspace{15mm}\\
   \end{array} \right. $$

3a) Let $n \ge 16$, then
$\dim H^2({\mathcal V}_n, {\mathcal V}_n)=10$,
more precisely:
$$\dim H^2_{(\mu)}({\mathcal V}_n, {\mathcal V}_n)
= \left\{\begin{array}{r}
   2, \hspace{25mm} \mu={-}2; \hspace{35mm}\\
   1, \quad \mu={-}4,{-}3,{-}1,n{-}11, n{-}10, n{-}9, n{-}8, n{-}7;
   \hspace{0mm}\\
   0, \hspace{26mm} otherwise. \hspace{30mm}\\
   \end{array} \right. $$

3b) Let $12 \le n \le 15$, then $\dim H^2({\mathcal V}_n, {\mathcal V}_n)=11$,
more precisely:
$$\dim H^2_{(\mu)}({\mathcal V}_n, {\mathcal V}_n)
= \left\{\begin{array}{r}
   2, \hspace{25mm} \mu={-}2; \hspace{35mm}\\
   1, \quad \mu={-}4,{-}3,{-}1,1,n{-}11, n{-}10, n{-}9, n{-}8, n{-}7;
   \hspace{0mm}\\
   0, \hspace{26mm} otherwise. \hspace{30mm}\\
   \end{array} \right. $$
\end{theorem}

\begin{proof}

For the proof we will use the spectral sequence considered above.
Let us recall some results from ~\cite{Mill2}, namely:

1) the basis of $H^1({\mathcal V}_n)$ 
consists of two classes $[e^1]$ and $[e^2]$.

2) the space $H^2({\mathcal V}_n)$ is $3$-dimensional
and generated
by classes 
$$g_5=[e^2 {\wedge} e^3],\;
g_7=[e^2 {\wedge} e^5 - 3 e^3 {\wedge} e^4], \;
[\Omega_{n{+}1}]{=}\frac{1}{2}[\sum\limits_
{i{+}j{=}n{+}1}(j{-}i)
e^i{\wedge}e^j]$$ of weights $5$, $7$, $n{+}1$ respectively.

3) $H^3({\mathcal V}_n)$ is $5$-dimensional and generated
by elements $g_{12}$, $g_{15}$ and
\begin{equation}
\begin{split}
[e^2{\wedge}\Omega_{n{+}1}],\:
[e^2{\wedge}\Omega_{n{+}2}-n e^3{\wedge}\Omega_{n{+}1}],\\
[e^2{\wedge}\Omega_{n{+}3}
-(n{+}1)e^3{\wedge}\Omega_{n{+}2}+
\frac{n(n{+}1)}{2}e^4{\wedge}\Omega_{n{+}1}]
\end{split}
\end{equation}
of weights $12$, $15$, $n{+}3$, $n{+}4$, $n{+}5$
respectively. Where
$$\Omega_{n{+}2}{=}\frac{1}{2}\sum\limits_
{{\scriptsize \begin{matrix}
i{+}j{=}n{+}2 \\
i,j>1
\end{matrix}}}
(j{-}i)e^i{\wedge}e^j, \quad
\Omega_{n{+}3}{=}\frac{1}{2}\sum\limits_
{{\scriptsize \begin{matrix}
i{+}j{=}n{+}3 \\
i,j>2
\end{matrix}}}
(j{-}i)e^i{\wedge}e^j$$
are projections of non-existing $de^{n{+}2}$ and $de^{n{+}3}$ to
$\Lambda^2({\mathcal V}_n)$. We have also
\begin{equation}
\begin{split}
d\Omega_{n{+}2}{=}d(de^{n{+}2}{-}ne^1{\wedge}e^{n{+}1}){=}
ne^1{\wedge}\Omega_{n{+}1},\\
d\Omega_{n{+}3}{=}d(de^{n{+}3}{-}(n{+}1)e^1{\wedge}e^{n{+}2}
{-}(n{-}1)e^2{\wedge}e^n){=}
(n{+}1)e^1{\wedge}\Omega_{n{+}2}{+}(n{-}1)e^2{\wedge}\Omega_{n{+}1}.
\end{split}
\end{equation}

\begin{remark}
The classes $g_{12}$, $g_{15}$ are generators of $H^3(L_1)$ of weights 
$12$, $15$ respectively, they arise in 
Gontcharova's theorem (see ~\cite{Fu} for details), as well
as $g_5$, $g_7$ in $H^2(L_1)$ and $[e^1]$, $[e^2]$ in $H^1(L_1)$.
We will need the following formula for $g_{12}$:
$$g_{12}=
[2e^2{\wedge} e^3{\wedge} e^7{-}
5e^2{\wedge} e^4{\wedge} e^6{+}20e^3{\wedge} e^4{\wedge} e^5].$$
\end{remark}

\begin{proposition}
\label{spdiff} Let $n\ge 16$,
then we have  only a finite number of non-trivial differentials
$d_r^{j,{-}j}, d_r^{j,1{-}j}, d_r^{j,2{-}j}$ of the spectral sequence
$E_r$. They are:
\begin{equation}
\label{diffs}
\begin{split}
d_1^{j,{-}j}:e_j \longmapsto (1{-}j)e_{j{+}1}{\otimes} [e^1],\:j \ne 1;\\
d_2^{1,1}:e_1 \longmapsto e_3{\otimes} [e^2], \\
d_2^{j,2{-}j}:
e_j {\otimes} [\Omega_{n{+}1}] \longmapsto
\left(j{-}2-\frac{(n{-}1)j(j{-}1)}{(n{+}1)n}\right)e_{j{+}2} {\otimes}
[e^2{\wedge}\Omega_{n{+}1}];\\
d_3^{j,1{-}j}:
e_j {\otimes} [e^2] \longmapsto
{-}\frac{1}{6}(j{-}3)(j^2{-}3j{+}8)e_{j{+}3} {\otimes}
[e^2{\wedge}e^3], j \ne 3;\\
d_4^{2,{-}1}:
e_2 {\otimes} [e^1] \longmapsto
{-}\frac{4}{3}e_6 {\otimes}
[e^2{\wedge}e^3];\\
d_5^{j{,}2{-}j}{:}
e_j {\otimes} [e^2 {\wedge} e^5 {-} 3 e^3 {\wedge} e^4]
\mapsto
{-}\frac{1}{5544}(j{-}8)(j^2{-}4j{+}27)(j^2{-}13j{+}48)
e_{j{+}5} {\otimes} g_{12}, j {\ne} 8;\\
d_8^{8,{-}6}:
{\rm Span}(e_8 {\otimes} [e^2 {\wedge} e^5 {-} 3 e^3 {\wedge} e^4])
\to
{\rm Span}(e_{16} {\otimes} g_{15}).\\
\end{split}
\end{equation}
If $12 \le n \le 15$ then only the differential $d_8^{8,{-}6}$ becomes trivial.

\end{proposition}

\begin{corollary}
\label{maincorl}
The following classes in $E_1^{p,q}, p{+}q{=}0,1,2$ survive to
$E_{\infty}$ for $n \ge 16$:
\begin{equation}
\begin{split}
e_n; \;
e_1 {\otimes} [e^1], e_n {\otimes} [e^2], e_{n{-}1} {\otimes} [e^2],
e_{n{-}2} {\otimes} [e^2];\\
e_{n{-}1} \otimes [\Omega_{n{+}1}], \qquad e_n \otimes [\Omega_{n{+}1}], \\
e_{i} \otimes [e^2{\wedge}e^3], \:i=1,2,3; \\
e_{j} \otimes [e^2 {\wedge} e^5 - 3 e^3 {\wedge} e^4], \:
j=n{-}4, n{-}3, n{-}2, n{-}1, n.
\end{split}
\end{equation}
If $12 \le n \le 15$ then one have to add  $e_{8} \otimes [e^2 {\wedge} e^5 - 3 e^3 {\wedge} e^4]$ to the list above.
\end{corollary}
\begin{remark}
1) the element $e_n$ spans the center
$Z({\mathcal V}_n)=H^0({\mathcal V}_n, {\mathcal V}_n)$;

2) the class $e_1 {\otimes} [e^1]$ corresponds to the inner derivation
$\sum_{i{=}1} ie_i {\otimes} e^i=ad(e_0)$
in the solvable Lie algebra $L_0/L_{n{+}1}$ restricted to the
nilpotent ideal ${\mathcal V}_n$.  
\end{remark}

Hence the proof of our theorem follows from Corollary~\ref{maincorl}.
It's time to say that the Proposition \ref{spdiff}
it is finite-dimensional version of the following theorem

\begin{theorem}[A.~Fialowski, ~\cite{Fial}, ~\cite{FiFu}]
$\dim H^2(L_1, L_1)=3$,
more precisely:
$$\dim H^2_{(\mu)}(L_1, L_1)
= \left\{\begin{array}{r}
   1, \quad \mu={-}2,{-}3,{-}4; \\
   0, \quad otherwise. \\
   \end{array} \right. $$
\end{theorem}

The differentials $d_1^{j,{-}j}$,$d_2^{1,1}$, $d_3^{j,1{-}j}$,
$d_4^{2,{-}1}$,
$d_5^{j,2{-}j}$, $d_8^{8,{-}6}$ of Proposition \ref{spdiff} came
from the corresponding spectral sequence for $H^*(L_1,L_1)$.
Their non-triviality follows from the more
general result of B.~Feigin, D.~Fuchs ~\cite{FeFu}. But it is
a very complicate task to follow all details of the proof in ~\cite{FeFu}.
So it appears to be usefull
to calculate $d_1^{j,{-}j}$, $d_3^{j,1{-}j}$, $d_5^{j,2{-}j}$
explicitly. From the other hand explicit formulae will give us
a possibility to write down structure relations of all deformations
of ${\mathcal V}_n$.

\begin{remark}
In the infinite-dimensional case only classes
$e_{i} \otimes [e^2{\wedge}e^3], \:i=1,2,3$
survive to $E_{\infty}$ and they correspond to generators
in $H^2(L_1,L_1)$ of weights $\mu=i-(2{+}3)={-}4,{-}3,{-}2$
that were found by A.~Fialowski.
\end{remark}

The proof of Proposition \ref{spdiff} consists of direct calculations.

First of all by definiton of 
$d:C^0(\mathfrak{g},\mathfrak{g}){=}\mathfrak{g} \to
C^1(\mathfrak{g},\mathfrak{g}){=}{\rm Hom}(\mathfrak{g},\mathfrak{g}){=}
\mathfrak{g} {\otimes} {\mathfrak g}^*$ we have
$$d(e_j) =
(j{-}1)e_{j{+}1} \otimes e^1+(j{-}2)e_{j{+}2} \otimes e^2+
(j{-}3)e_{j{+}3} \otimes e^3+...$$
Hence $$d_1^{j,{-}j}(e_j) = e_{j{+}1} \otimes (j{-}1)[e^1], \quad j \ne 1.$$

Now we go to the differential $d_2$.
$$d(e_j \otimes \Omega_{n{+}1})=
e_{j{+}1} \otimes (j{-}1)e^1{\wedge}\Omega_{n{+}1}+...$$
where dots stand instead of terms of higher filtration.
As we know $e^1{\wedge}\Omega_n=\frac{1}{n}d\Omega_{n{+}2}$ and
we now can take new representative $e_j \otimes \Omega_n+
\frac{j{-}1}{n}e_{j{+}1} \otimes \Omega_n$ such that:
$$d(e_j {\otimes} \Omega_{n{+}1}+\frac{j{-}1}{n}e_{j{+}1}
{\otimes} \Omega_{n{+}2})
=e_{j{+}2} {\otimes} \left( \frac{j(j{-}1)}{n}e^1{\wedge}\Omega_{n{+}2}
+(j{-}2)e^2 {\wedge} \Omega_{n{+}1} \right){+}\dots$$

But $d\Omega_{n{+}3}=(n{+}1)e^1{\wedge}\Omega_{n{+}2}+
(n{-}1)e^2{\wedge}\Omega_{n{+}1}$ and therefore
$\frac{j(j{-}1)}{n}e^1{\wedge}\Omega_{n{+}2}
+(j{-}2)e^2 {\wedge} \Omega_{n{+}1}$ is cohomologous to
$\left(j{-}2-\frac{(n{-}1)j(j{-}1)}{(n{+}1)n}\right)
e^2{\wedge}\Omega_{n{+}1}$. Hence we conclude

$$d_2(e_j {\otimes} [\Omega_{n{+}1}])
=e_{j{+}2} {\otimes} \left(j{-}2-\frac{(n{-}1)j(j{-}1)}{(n{+}1)n}\right)
[e^2 {\wedge} \Omega_{n{+}1}].$$
 
Now the most complicated case: computation of
$d_2(e_{j} \otimes [e^2 {\wedge} e^5 - 3 e^3 {\wedge} e^4])$.
We have to find $\xi_1, \xi_2, \xi_3, \xi_4$ such that
\begin{equation}
\label{d5}
d(e_{j} {\otimes} (e^2 {\wedge} e^5 {-} 3 e^3 {\wedge} e^4)+
\sum_{p{=}1}^4 e_{j{+}p} {\otimes} \xi_p)= 
e_{j{+}5} {\otimes} \xi +\dots
\end{equation}
for some $\xi \in C^3({\mathcal V}_n, {\mathcal V}_n)$.
We recall the notation $g_7=e^2 {\wedge} e^5 {-} 3 e^3 {\wedge} e^4$.
We have the following system of equations on $\xi_1, \xi_2, \xi_3, \xi_4$:
\begin{equation}
\begin{split}
d\xi_1=(j{-}1)e^1{\wedge}g_7;\\
d\xi_2=je^1{\wedge}\xi_1+(j{-}2)e^2{\wedge}g_7;\\
d\xi_3=(j{+}1)e^1{\wedge}\xi_2+(j{-}2)e^2{\wedge}\xi_1+
(j{-}3)e^3{\wedge}g_7;\\
d\xi_4=(j{+}2)e^1{\wedge}\xi_3+je^2{\wedge}\xi_2+
(j{-}2)e^2{\wedge}\xi_1+(j{-}4)e^4{\wedge}g_7.\\
\end{split}
\end{equation}

Taking $\xi_1, \xi_2, \xi_3, \xi_4$ as homogeneous $2$-forms
of weights $8,9,10,11$ one can remark that
the right parts of these equations are exact forms because
$H^3_{(p)}({\mathcal V}_n)=0, p \le 11$.
$\xi_1, \xi_2, \xi_3, \xi_4$ are defined uniquely by the condition
$\xi_p \in \Lambda^2(e^2,\dots,e^n), \: p=1,2,3,4$.

The answer is:
\begin{equation}
\begin{split}
\xi_1=\frac{j{-}1}{2}(e^2{\wedge}e^6-2e^3{\wedge}e^5);\\
\xi_2=P_1(j)e^2{\wedge}e^7+P_2(j)e^3{\wedge}e^6+
P_3(j)e^4{\wedge}e^5;\\
\xi_3=Q_1(j)e^2{\wedge}e^8+Q_2(j)e^3{\wedge}e^7+
Q_3(j)e^4{\wedge}e^6;\\
\xi_4=Z_1(j)e^2{\wedge}e^9+Z_2(j)e^3{\wedge}e^8+
Z_3(j)e^4{\wedge}e^7+Z_4(j)e^5{\wedge}e^6,\\
\end{split}
\end{equation}
where polinomials $P_i(j), Q_j(j), Z_l(j)$ are defined by
\begin{equation}
\begin{split}
P_1(j){=}\frac{\left(5\binom{j}{2}{+}3\binom{j}{1}{-}6\right)}{21};\:
P_2(j){=}{-}\frac{\left(4\binom{j}{2}{+}15\binom{j}{1}{-}30\right)}{21};\:
P_3(j){=}{-}\frac{\left(13\binom{j}{2}{-}30\binom{j}{1}{+}60\right)}{21};\\
Q_1(j){=}\frac{\left(3\binom{j{+}1}{3}+4\binom{j{+}1}{2}{-}
4\binom{j{+}1}{1}\right)}{28};\:
Q_2(j){=}\frac{\left(\binom{j{+}1}{3}{-}8\binom{j{+}1}{2}{+}
8\binom{j{+}1}{1}\right)}{14};\\
Q_3(j){=}\frac{\left({-}13\binom{j{+}1}{3}{+}20\binom{j{+}1}{2}{-}
20\binom{j{+}1}{1}\right)}{28};\\
Z_1(j){=}\frac{1}{22}\binom{j{+}2}{4}{+}\frac{23}{231}\binom{j{+}2}{3}{-}
\frac{7}{198}\binom{j{+}2}{2}{-}\frac{59}{693}\binom{j{+}2}{1}{+}
\frac{37}{231};\\
Z_2(j){=}\frac{17}{154}\binom{j{+}2}{4}{-}\frac{62}{231}\binom{j{+}2}{3}{-}
\frac{53}{1386}\binom{j{+}2}{2}{-}\frac{59}{99}\binom{j{+}2}{1}{-}
\frac{37}{33};\\
Z_3(j){=}{-}\frac{29}{154}\binom{j{+}2}{4}{-}\frac{4}{77}\binom{j{+}2}{3}{+}
\frac{317}{462}\binom{j{+}2}{2}{-}\frac{59}{33}\binom{j{+}2}{1}{+}
\frac{37}{11};\\
Z_4(j){=}{-}\frac{47}{154}\binom{j{+}2}{4}{+}\frac{185}{231}
\binom{j{+}2}{3}{-}
\frac{2245}{1386}\binom{j{+}2}{2}{+}\frac{295}{99}\binom{j{+}2}{1}{-}
\frac{185}{33};
\end{split}
\end{equation}
Now one can calculate $\xi$ from equation (\ref{d5}):
$$
\xi=(j{+}3)e^1{\wedge}\xi_4+(j{+}1)e^2{\wedge}\xi_3+
(j{-}1)e^3{\wedge}\xi_2+(j{-}3)e^4{\wedge}\xi_1+(j{-}5)e^5{\wedge}g_7.
$$
The space $H^3_{(12)}({\mathcal V_n})$ is one-dimensional and
we have for the cohomology class $[\xi]$:
$$[\xi]=
{-}\frac{1}{5544}(j{-}8)(j^2{-}4j{+}27)(j^2{-}13j{+}48)
\left[2e^2{\wedge} e^3{\wedge} e^7{-}
5e^2{\wedge} e^4{\wedge} e^6{+}20e^3{\wedge} e^4{\wedge} e^5\right].$$

In the same way one can remark that
\begin{equation}
\begin{split}
d(e_{j} {\otimes} e^2 {+}
e_{j{+}1} {\otimes} (j{-}1)e^3{+}
e_{j{+}2} {\otimes} \frac{j(j{-}1)}{2}e^4)
=\\
=e_{j{+}3} {\otimes}\left( \frac{(j{+}1)j(j{-}1)}{2}e^1{\wedge}e^4+
\left((j{-}1)^2{-}(j{-}3)\right)e^2{\wedge}e^3\right)
+\dots
\end{split}
\end{equation}
As $3e^1{\wedge}e^4$ is cohomologous to $e^2{\wedge}e^3$ it follows that
$$
d_3(e_{j} {\otimes} [e^2])
=e_{j{+}3} {\otimes}
\left((j{-}1)^2{-}(j{-}3)-\frac{(j{+}1)j(j{-}1)}{6}\right)[e^2{\wedge}e^3]
$$
\end{proof}

Now we choose the following basic cocycles $\psi_{k,i}$
in $\oplus_{\mu>0}
H^2_{(\mu)}({\mathcal V}_n, {\mathcal V}_n), \:n \ge 16$:

$$
\begin{tabular}{|c|c|}
\hline
&\\[-10pt]
& basic cocycles\\
&\\[-10pt]
\hline
&\\[-10pt]
$H^2_{(n{-}7)}({\mathcal V}_n, {\mathcal V}_n)$ & \begin{tabular}{c}
$\psi_{n,7}=e_n\otimes(e^2 {\wedge} e^5 - 3 e^3 {\wedge} e^4)$\\
\end{tabular}\\
&\\[-10pt]
\hline
&\\[-10pt]
$H^2_{(n{-}8)}({\mathcal V}_n, {\mathcal V}_n)$ &
$\psi_{n,8}=e_{n{-}1}\otimes (e^2 {\wedge} e^5 - 3 e^3 {\wedge} e^4)+
e_{n}\otimes \frac{n{-}2}{2} (e^2 {\wedge} e^6 - 2 e^3 {\wedge} e^5)$;
\\
&\\[-10pt]
\hline
&\\[-10pt]
$H^2_{(n{-}9)}({\mathcal V}_n, {\mathcal V}_n)$ & \begin{tabular}{c}
$\psi_{n,9}=e_{n{-}2}\otimes (e^2 {\wedge} e^5 - 3 e^3 {\wedge} e^4)+
e_{n{-}1}\otimes \frac{n{-}3}{2} (e^2 {\wedge} e^6 - 2 e^3 {\wedge} e^5)+$\\
$+ \: e_n\otimes(P_1(n{-}2)e^2 {\wedge} e^7+P_2(n{-}2)e^3 {\wedge} e^6
+P_3(n{-}2)e^4 {\wedge} e^5)$\\
\end{tabular}\\
&\\[-10pt]
\hline
&\\[-10pt]
$H^2_{(n{-}10)}({\mathcal V}_n, {\mathcal V}_n)$ & \begin{tabular}{c}
$\psi_{n,10}=e_{n{-}3}\otimes (e^2 {\wedge} e^5 - 3 e^3 {\wedge} e^4)+
e_{n{-}2}\otimes \frac{n{-}4}{2} (e^2 {\wedge} e^6 - 2 e^3 {\wedge} e^5)+$\\
$+ \: e_{n{-}1}\otimes(P_1(n{-}3)e^2 {\wedge} e^7+P_2(n{-}3)e^3 {\wedge} e^6
+P_3(n{-}3)e^4 {\wedge} e^5)+$\\
$+ \: e_{n}\otimes(Q_1(n{-}3)e^2 {\wedge} e^8+Q_2(n{-}3)e^3 {\wedge} e^7
+Q_3(n{-}3)e^4 {\wedge} e^6)$\\
\end{tabular}\\
&\\[-10pt]
\hline
&\\[-10pt]
$H^2_{(n{-}11)}({\mathcal V}_n, {\mathcal V}_n)$ & \begin{tabular}{c}
$\psi_{n,11}=e_{n{-}4}\otimes (e^2 {\wedge} e^5 - 3 e^3 {\wedge} e^4)+
e_{n{-}3}\otimes \frac{n{-}5}{2} (e^2 {\wedge} e^6 - 2 e^3 {\wedge} e^5)+$\\
$+ \: e_{n{-}2}\otimes(P_1(n{-}4)e^2 {\wedge} e^7+P_2(n{-}4)e^3 {\wedge} e^6
+P_3(n{-}4)e^4 {\wedge} e^5)+$\\
$+ \: e_{n{-}1}\otimes(Q_1(n{-}4)e^2 {\wedge} e^8+Q_2(n{-}4)e^3 {\wedge} e^7
+Q_3(n{-}4)e^4 {\wedge} e^6)+$\\
$e_{n}{\otimes}(Z_1(n{-}4)e^2 {\wedge} e^9{+}Z_2(n{-}4)e^3 {\wedge} e^8
{+}Z_3(n{-}4)e^4 {\wedge} e^7{+}Z_4(n{-}4)e^5 {\wedge} e^6)$\\
\end{tabular}\\
\hline
\end{tabular}$$

\section{Moduli space of ${\mathbb Z}_{{>}0}$-filtered deformations.}
In this section we classify up to an isomorphism the Lie algebras over ${\mathbb K}$ defined
by the basis $e_1,\dots, e_n, \: n\ge 16$ and 
commutating relations of the following form:

\begin{equation}
\label{defgonch}
[e_i, e_j]= (j-i)e_{i{+}j}+\sum \limits_{l{=}1}^{n{-}i{-}j}c_{ij}^{l}e_{i{+}j{+}l}.   
\end{equation}

A Lie algebra $\mathfrak g$ with the commutating relations 
(\ref{defgonch}) is a 
${\mathbb Z}_{{>}0}$-filtered deformation of 
the ${\mathbb Z}_{{>}0}$-graded Lie algebra ${\mathcal V}_n$, i.e.
$\mathfrak g=({\mathcal V}_n, [,]+\Psi)$,
where    
\begin{equation}
\label{defcond}
\begin{split}
\Psi=\Psi_1+\Psi_2+\dots+\Psi_{n{-}3}, \quad \Psi_l \in C^2_{(l)}({\mathcal V}_n, {\mathcal V}_n),\\
\Psi_l(e_i,e_j)=\left\{ \begin{array}{l} c_{ij}^{l} e_{i{+}j{+}l}, \;\; i{+}j \le n{-}l,\\ 
0, \quad {\rm otherwise}; \end{array} \right.
\quad l=1,2,\dots,n{-}3.\\
\end{split}
\end{equation}
satisfying to the system (\ref{sysdef}) of deformation 
equations: 
$$
d\Psi_1=0, \;
d\Psi_{2}+\frac{1}{2}[\Psi_{1},\Psi_{1}]=0, \;
\dots ,
d\Psi_{n{-}6}+\frac{1}{2}\sum \limits_{i{+}j{=}n{-}6}[\Psi_i,\Psi_j]=0.
$$

\begin{lemma}
\label{isom_lemma}
Let $\Psi, \tilde \Psi$ be two ${\mathbb Z}_{{>}0}$-filtered deformations of 
${\mathcal V}_n, n \ge 5$ and
$$\varphi : \mathfrak{g}=({\mathcal V}_n,[,]+\Psi ) \to 
\tilde {\mathfrak{g}} =({\mathcal V}_n,[,]+\tilde \Psi)$$ 
is a Lie algebra isomorphism.

Then 
\begin{equation}
\begin{split}
1) \; \varphi=\varphi_0+\varphi_1+\varphi_2+\dots+\varphi_{n{-}1}, \quad 
\varphi_j \in C^1_{(j)}({\mathcal V}_n, {\mathcal V}_n), \; j=0,1,2,\dots,n{-}1, \\ 
2) \; \varphi_0(e_i)=\alpha^ie_i, \: \alpha \in {\mathbb K}^*, \quad
i=1,2,\dots,n. \hspace{5.15cm}
\end{split}
\end{equation}
\end{lemma}

\begin{proof}
We will study the dual situation. The mapping 
$$\varphi^*: {\mathfrak g}^*=({\mathcal V}_n^*,d_{\Psi}=d+\Psi^*) \to 
\tilde {\mathfrak g}^*=({\mathcal V}_n^*,d_{\tilde \Psi}=d+{\tilde \Psi}^*)$$
is an isomorphism of $d$-algebras.

For the dual basis $e^1,e^2,\dots, e^n$ of $({\mathcal V}_n^*, d_{\Psi})$ we have the following
structure relations:  
$$
d_{\Psi}e^k= \frac{1}{2}\sum_{i{+}j{=}k}
(j-i)e^i \wedge e^j+\frac{1}{2}\sum_{m{+}p{<}k}
c_{mp}^k e^m \wedge e^p, \; k=3,\dots,n.
$$

Let us write down some of them. 
\begin{equation}
\begin{split}
d_{\Psi}e^1=d_{\Psi}e^2=0, \; d_{\Psi}e^3=e^1 \wedge e^2, \;
d_{\Psi}e^4=2e^1 \wedge e^3+c_{12}^1e^1 \wedge e^2, \\  
d_{\Psi}e^5=3e^1 \wedge e^4+e^2 \wedge e^3+c_{13}^1 e^1 \wedge e^3+ c_{12}^2e^1 \wedge e^2,\dots
\end{split}
\end{equation}

The dual mapping $\varphi^* : {\tilde {\mathfrak g}}^* \to  {\mathfrak g}^*$
is the isomorphism of $d$-algebras, i.e. $$ d_{\Psi} \varphi^* = \varphi^* d_{\tilde \Psi}$$

1)  $ \; d_{\Psi}\varphi^* e^1=d_{\Psi}\varphi^* e^2=0$, thus 
$$
\varphi^* e^1=\alpha_{11}e^1+\alpha_{21}e^2;\quad
\varphi^* e^2=\alpha_{12}e^1+\alpha_{22}e^2.
$$

2)$\; d_{\Psi}\varphi^* e^3=\varphi^* e^1 \wedge \varphi^* e^2=
(\alpha_{11}\alpha_{22}-\alpha_{21}\alpha_{12})e^1 \wedge e^2$,
and we have
$$
\varphi^* e^3=(\alpha_{11}\alpha_{22}-\alpha_{21}\alpha_{12})e^3+\alpha_{13}e^1+\alpha_{23}e^2.
$$

3) The form $\varphi^* d_{\tilde \Psi}e^4$ is cohomologous to 
$2\alpha_{21}(\alpha_{11}\alpha_{22}-\alpha_{21}\alpha_{12})e^2 \wedge e^3$
and it is exact iff $\alpha_{21}=0$, thus
$$
\varphi^* e^4=\alpha_{11}^2\alpha_{22} e^4 +(\alpha_{11}\alpha_{23}+
{\tilde c}_{12}^1 \alpha_{11}\alpha_{22}-
c^1_{12}\alpha_{11}^2\alpha_{22})e^3+
\alpha_{14}e^1+\alpha_{24}e^2.
$$

4) $\varphi^* d_{\tilde \Psi}e^5 
\; \sim \; 3\alpha_{22}\alpha_{11}^3e^1 \wedge e^4+\alpha_{22}^2\alpha_{11}e^2 \wedge e^3$
and the last one is exact iff $\alpha_{12}^2=\alpha_{22}$, and for the moment we have
\begin{equation}
\begin{split}
\varphi^* e^1=\alpha_{11}e^1, \quad \varphi^* e^2=\alpha_{11}^2e^2+\alpha_{12}e^1,\\
\varphi^* e^3=\alpha_{11}^3e^3+\alpha_{13}e^1+\alpha_{23}e^2,\quad\\
\varphi^* e^4=\alpha_{11}^4e^4+\dots, \quad
\varphi^* e^5=\alpha_{11}^5e^5+\dots
\end{split}
\end{equation}

Going on and using an obvious inductive assumption we have for the operator
$\varphi^*$:
\begin{equation}
\label{triang_m}
\varphi^*e^i=\alpha^i e^i+\sum_{l<i}^n \alpha_{li}e^l, \; i=1,\dots, n,
\end{equation}
for some $\alpha \ne 0, \alpha_{li} \in {\mathbb K}$.
\end{proof}

From the other hand it is evident that changing the canonical basis of an arbitrary Lie algebra 
$\mathfrak g$  
of the type (\ref{defgonch}) by an operator $\varphi$ with property (\ref{triang_m})
we will get again the commutating relations of the same type.
 
\begin{corollary}
The matrix Lie group $G_n$ of lower-triangular matrixes $\varphi$ of the following type 
$$\varphi=
 \left( \begin{array}{cccc}
\alpha & 0& \ldots & 0\\
a_{21} & \alpha^2 & \ldots & 0\\
\vdots & \vdots & \ddots & \vdots\\
a_{n1} & a_{nn{-}1} & \ldots & \alpha^n
\end{array} \right), \;a_{ij},\: \alpha \in {\mathbb K}, \:\alpha \ne 0,
$$
acts on the set $V_n$ of
${\mathbb Z}_{{>}0}$-filtered deformations of ${\mathcal V}_n$
as the group of changes of canonical basis:
\begin{equation}
\label{G_n-action}
(\varphi \star \Psi)(x,y)=\varphi^{{-}1}([\varphi x, \varphi y]+\Psi(\varphi x, \varphi y))-[x,y],
\; \forall x,y \in {\mathcal V_n}, \; \varphi \in G_n. 
\end{equation}
Two ${\mathbb Z}_{{>}0}$-filtered deformations
$({\mathcal V}_n,[,]+\Psi )$ and $({\mathcal V}_n,[,]+\tilde \Psi )$
are isomorphic as Lie algebras if and only if they are in the same
orbit $O_{\Psi}$ of the $G_n$-action. 
\end{corollary}

Hence we have proved the following
\begin{theorem} Let $n \ge 5$, then 
there is a one-to-one correspondance between the orbit space $O(G_n,V_n)$ 
of the action  $G_n$ on the set $V_n$  of ${\mathbb Z}_{{>}0}$-filtered deformations of ${\mathcal V}_n$
and the moduli space ${\mathcal M}_n$  
(the set of isomorphism classes of the ${\mathbb Z}_{{>}0}$-filtered deformations of ${\mathcal V}_n$). 
\end{theorem}

\begin{proposition} The matrix Lie group $G_n$ is
the semi-direct product ${\mathbb K}^* \ltimes UT_n$  where $UT_n$ denotes the group of 
unitriangular matrixes:
$$G_n  = {\mathbb K}^* \ltimes UT_n=
{\mathbb K}^* \ltimes 
\left\{ \left( \begin{array}{cccc}
1 & 0& \ldots & 0\\
* & 1 & \ldots & 0\\
\vdots & \vdots & \ddots & \vdots\\
* & * & \ldots & 1
\end{array} \right) \right\}, \;
{\mathbb K}^*\cong 
\left\{
 \left( \begin{array}{cccc}
\alpha & 0& \ldots & 0\\
0 & \alpha^2 & \ldots & 0\\
\vdots & \vdots & \ddots & \vdots\\
0 & 0 & \ldots & \alpha^n
\end{array} \right) \right\}.
$$
\end{proposition}

\begin{remark}[see also section \ref{Nijenhuis-Richardson}]
If $\Psi=\Psi_1+\Psi_2+\dots+\Psi_{n{-}3}$ is a solution of the system 
(\ref{sysdef}) of deformation equations then 
$\Psi=t\Psi_1+t^2\Psi_2+\dots+t^{n{-}3}\Psi_{n{-}3},\; \forall t \in {\mathbb K}$
also satisfies to the system 
(\ref{sysdef}). It follows that 
the space $V_n$ can be retracted over itself to ${\mathcal V}_n$.   
From another hand  one can define ${\mathbb K}^*$-action on $V_n$:
$$\rho_n( \alpha)
(\Psi_1,\Psi_2,\dots,\Psi_{n{-}3})=(\alpha\Psi_1,\alpha^2\Psi_2,\dots,\alpha^{n{-}3}\Psi_{n{-}3}), \; \alpha \in 
{\mathbb K}^*.$$ Evidently this action coincides with the action on $V_n$ of the subgroup ${\mathbb K}^*$ 
of diagonal matrixes in $G_n$. 
\end{remark}

\begin{proposition}
There are the following bijections:
$$O_{\rho_n}({\mathbb K}^*,O(UT_n,V_n)) \rightarrow O(G_n/UT_n,O(UT_n,V_n)) \rightarrow O(G_n,V_n).$$ 
\end{proposition} 

\begin{proposition} 
Let $\Psi$ be a ${\mathbb Z}_{{>}0}$-filtered deformation of ${\mathcal V}_n, n \ge 12$.
Then there exists an element $\tilde \Psi$ 
in the $UT_n$-orbit $O_{\Psi}$ of $\Psi$ such that
$$\tilde \Psi_1=\dots=\tilde \Psi_{n{-}12}=0.$$
\end{proposition}

For an arbitrary ${\mathbb Z}_{{>}0}$-filtered deformation $\Psi$ the first equation is
$d\Psi_1=0$. 
We recall now  
$$H^2_{(i)}({\mathcal V}_n,{\mathcal V}_n)=0,\:i \le n{-}12$$ 
and hence $\Psi_1=d\varphi_1$ for some $\varphi_1$  of
$C^1_{(1)}({\mathcal V}_n, {\mathcal V}_n)$. Acting by $g=id+\varphi_1 \in UT_n$
we get $\tilde \Psi=g \star \Psi$ such that $\tilde \Psi_1=0$. 
Now  the second equation of the system for this new element $\tilde \Psi$.  
will be $d\tilde \Psi_2=[\tilde\Psi_1,\tilde \Psi_1]=0$. We act on $\tilde \Psi$ by $g=id+\varphi_2$,
where $\tilde \Psi_2=d\varphi_2$ and continue the procedure step by step.   

Now we suppose to be constructed an element $\tilde \Psi$ such that 
$\tilde \Psi_1=\dots=\tilde \Psi_{n{-}12}=0$ satisfying  
$$d\tilde \Psi_{n{-}11}=0.$$

We recall that in the section \ref{Computations} we found the basic cocycles $\psi_{n,i}$, such that 
$$Span([\psi_{n,i}])=H^2_{(n{-}i)}({\mathcal V}_n,{\mathcal V}_n),\;
i=7,\dots,11.$$

Hence $\tilde \Psi_{n{-}11}=x_1\psi_{n,11}+d\varphi_{n{-}11}, x_1 \in {\mathbb K}$, 
and $\varphi_{n{-}11} \in C^1_{(n{-}11)}({\mathcal V}_n,{\mathcal V}_n))$.  
Then  acting on $\tilde \Psi$  by $id+\varphi_{n{-}11}$
we get again a new element in the orbit $O_{\Psi}$ 
(we keep the same notation $\tilde \Psi$ for it ) with the property
$\tilde \Psi_{n{-}11}=x_1 \psi_{n,11}$.  

\begin{proposition}
Let  $n\ge 14$ then all Nijenhuis-Richardson products of basic cocycles $\psi_{n,i}$
are trivial elements in $C^3({\mathcal V}_n, {\mathcal V}_n)$:
$$[\psi_{n,i},\psi_{n,j}](x,y,z)=0 \: \; \forall x,y,z \in {\mathcal V}_n.$$
\end{proposition}

The proof follows from the 
two properties of $\psi_{n,l}$:

1) ${\rm Im}\psi_{n,l}={\rm Span}(e_n,\dots,e_{n{-}4});$

2) $\psi_{n,l}(x,y)=0$ if $x\wedge y \notin \Lambda^2(e_2,\dots,e_9)$.

Hence $\psi_{n,i}(\psi_{n,j}(x,y),z)=0, \: \forall x,y,z$
if $ n{-}4 > 9$.

\begin{proposition} 
Let $n \ge 14$. 
There is a one-to-one correspondance between the orbit space $O(UT_n, V_n)$ and
the $5$-dimensional vector space $\oplus_{i>0}H^2_{(i)}({\mathcal V}_n, {\mathcal V}_n)$.
\end{proposition}
This proposition follows from the previous one. 
Namely in an arbitrary $UT_n$-orbit  
$O_{\Psi}$ one can choose   
the unique representative
$\tilde \Psi$ such that
$$\tilde \Psi=x_1\psi_{n,11}+x_2\psi_{n,10}+
x_3\psi_{n,9}+x_4\psi_{n,8}+ x_5\psi_{n,7},$$
where $x_i \in {\mathbb K}, i=1,\dots,5$.
We will call $\tilde \Psi$ the canonical element of the orbit $O_{\Psi}$ and
the set $\{x_1,\dots,x_5\}$ is called the homogeneous coordinates of the orbit $O_{\Psi} \in O(UT_n,V_n)$. 

\begin{theorem}
\label{main}
Let $n \ge 16$. 
There is a one-to-one correspondance between the moduli space ${\mathcal M}_n=O(G_n, V_n)$ 
of ${\mathbb Z}_{{>}0}$-filtered deformations of  ${\mathcal V}_n$
and the orbit space $O_{\tilde \rho_n}({\mathbb K}^*, {\mathbb K}^5)$ where
the action $\tilde \rho_n$ of ${\mathbb K}^*$ on ${\mathbb K}^5$ is defined in coordinates $x^i$:
\begin{equation}
 \tilde \rho_n(\alpha) (x_1,x_2,\dots,x_5)=
(\alpha^{n{-}11}x_1,\alpha^{n{-}10}x_2,\dots,\alpha^{n{-}7}x_5),
 \; \alpha \in {\mathbb K}^*. 
\end{equation}
\end{theorem}
We have to verify only the formula for ${\mathbb K}^*$-action on $O(UT_n, V_n)$.
But 
\begin{equation}
\begin{split}
\tilde \rho_n(\alpha) (x_1,x_2,\dots,x_5)=\rho_n(\alpha) (x_1\psi_{n,11}+x_2\psi_{n,10}+\dots+x_5\psi_{n,17})=\\
=\alpha^{n{-}11}x_1\psi_{n,11}+\alpha^{n{-}10}x_2\psi_{n,10}+\dots+
\alpha^{n{-}7}x_5\psi_{n,7}.
\end{split}
\end{equation}

\section{Affine variety of ${\mathbb Z}_{{>}0}$-filtered deformations}
One can regard the set $V_n$ of ${\mathbb Z}_{{>}0}$-filtered deformation of 
${\mathcal V}_n$ as a affine variety in the affine space $\oplus_{i{>}0}C^2_{(i)}({\mathcal V}_n, {\mathcal V}_n)$. 
Namely one can define a mapping (polynomial in the coordinates $\{c_{ij}^k\}$):
$${\mathcal F}: \oplus_{i{>}0}C^2_{(i)}({\mathcal V}_n, {\mathcal V}_n) \rightarrow
\oplus_{i{>}0}C^3_{(i)}({\mathcal V}_n, {\mathcal V}_n),$$
where
$${\mathcal F}(\Psi_1,\Psi_2,\dots,\Psi_{n{-}3})=
(d\Psi_1,\:d\Psi_2+\frac{1}{2}[\Psi_1,\Psi_1],\:\dots,\:d\Psi_{n{-}6}+\frac{1}{2}\sum \limits_{i{+}j{=}n{-}6}[\Psi_i,\Psi_j]).$$
Then
$$V_n=\{\Psi \in \oplus_{i{>}0}C^2_{(i)}({\mathcal V}_n, {\mathcal V}_n)\; | \; \;{\mathcal F}(\Psi)=0\}.$$

\begin{theorem}
\begin{enumerate}
\item The variety $V_n$ has no singular points and 
$$\dim V_n= \left\{ \begin{array}{l} \frac{n(n{-}3)}{2}+3, \; n\ge 16;\\
\frac{n(n{-}3)}{2}+4, \; 12 \le n \le 15. \end{array} \right.$$  
\item The group $UT_n$ acts on $V_n$ with constant rank  and we have for 
the dimension of an arbitrary orbit $O_{\Psi}$
$$\dim O_{\Psi}= \frac{n(n{-}3)}{2}-2.$$  
\end{enumerate}
\end{theorem}

One can identify the Lie algebra of $UT_n$ with $\oplus_{{>}0}C^1_{(i)}({\mathcal V}_n, {\mathcal V}_n)$.
Let $\Psi$ be some fixed element of $V_n$, then 
decomposing 
$\varphi=\exp \alpha =1+\alpha+\alpha^2+ \dots$ in the formula (\ref{G_n-action}) and taking linear terms
with respect to $\alpha$ one can get the following formula for the differential $D$ of
the $UT_n$-action ($x,y \in {\mathcal V}_n$):
$$D_{\Psi}(\alpha)(x,y)=[\alpha(x),y]{+}[x,\alpha(y)]{-}\alpha([x,y])+\Psi(\alpha(x),y){+}\Psi(x,\alpha(y)){-}\alpha(\Psi(x,y)).
$$
Hence $$D_{\Psi}(\alpha)=-d(\alpha)+[\Psi,\alpha]=-d_{\Psi}(\alpha), \; \alpha \in 
\oplus_{{>}0}C^1_{(i)}({\mathcal V}_n, {\mathcal V}_n)$$ where 
$d_{\Psi} :\oplus_{{>}0}C^1_{(i)}({\mathcal V}_n, {\mathcal V}_n) \to \oplus_{{>}0}C^2_{(i)}({\mathcal V}_n, {\mathcal V}_n)$
is a deformed differential of the cochain complex of ${\mathcal V}_n$. But in the same time $d_{\Psi}$ defines the differential 
of the cochaine complex of the Lie algebra $\mathfrak g=({\mathcal V}_n,[,]+\Psi)$.

Analogously one can show that the differential $D{\mathcal F}_{\Psi}$ coincides with the differential
$d_{\Psi} :\oplus_{{>}0}C^2_{(i)}({\mathcal V}_n, {\mathcal V}_n) \to \oplus_{{>}0}C^3_{(i)}({\mathcal V}_n, {\mathcal V}_n)$
of the cochain complex with coefficients in the adjoint representation of $\mathfrak g=({\mathcal V}_n,[,]+\Psi)$. 

\begin{remark}
Let us denote by 
$d_i:C^*_{(i)}({\mathcal V}_n, {\mathcal V}_n) \to C^{*{+}1}_{(i)}({\mathcal V}_n, {\mathcal V}_n)$ 
the restriction of the differential $d$ of the cochain complex 
$(C^*({\mathcal V}_n, {\mathcal V}_n),d)$ to homogeneous components. 
Then the matrix of $d_{\Psi}$ for an arbitrary $\Psi$ is a block-triangular for some natural number $l$:
$$d_{\Psi}=
 \left( \begin{array}{cccc}
d_1 & *& \ldots & *\\
0 & d_2 & \ldots & *\\
\vdots & \vdots & \ddots & \vdots\\
0 & 0 & \ldots & d_{n{-}l}
\end{array} \right) 
$$
And we have the following estimates
$${\rm rank}\: d_{\Psi} \ge {\rm rank}\: d, \quad \dim {\rm Ker}\: d_{\Psi} \le \dim {\rm Ker}\: d.$$
\end{remark}  
In order to compute these ranks one may assume that $\Psi= \tilde \Psi$, where $\tilde \Psi \in {\mathbb K}^5$ is the canonical 
represantative of the orbit $O_{\Psi}$.

\begin{proposition}
\begin{equation}
\begin{split}
\dim {\rm Ker}\: \left( d_{\tilde \Psi}:\oplus_{{>}0}C^1_{(i)}({\mathcal V}_n, {\mathcal V}_n) \to 
\oplus_{{>}0}C^2_{(i)}({\mathcal V}_n, {\mathcal V}_n) \right)=\\
= \dim {\rm Ker}\: \left( d :\oplus_{{>}0}C^1_{(i)}({\mathcal V}_n, {\mathcal V}_n) \to 
\oplus_{{>}0}C^2_{(i)}({\mathcal V}_n, {\mathcal V}_n)\right)
= n+2.
\end{split}
\end{equation}
\end{proposition} 
One can easily verify that operators
\begin{equation}
\label{basis_yadra_stab}
\begin{split}
ad_{\tilde \Psi}(e_1),ad_{\tilde \Psi}(e_2),\dots, ad_{\tilde \Psi}(e_{n{-}1}), 
e_n{\otimes} e^2, e_{n{-}1} {\otimes} e^2+(n{-}2)e_n{\otimes }e^3,\\ e_{n{-}2} {\otimes} e^2+(n{-}3)e_n{\otimes} e^3+
\frac{(n{-}2)(n{-}3)}{2}e_n{\otimes} e^3 \hspace{2cm} 
\end{split}
\end{equation}
give the basis of ${\rm Ker}\: d_{\tilde \Psi}: 
\oplus_{{>}0}C^1_{(i)}({\mathcal V}_n, {\mathcal V}_n) \to \oplus_{{>}0}C^2_{(i)}({\mathcal V}_n, {\mathcal V}_n)$.

Hence $$\dim O_{\Psi}=\dim {\rm Im}\:\left( d_{\Psi}: 
\oplus_{{>}0}C^1_{(i)}({\mathcal V}_n, {\mathcal V}_n) \to \oplus_{{>}0}C^2_{(i)}({\mathcal V}_n, {\mathcal V}_n)\right)=
\frac{n(n{-}1)}{2}{-}n{-}2.$$
From the other hand the cocycles $\{\psi_{n,i}| i{=}7,{\dots},11\}$ span the tangent space to ${\mathbb K}^5$
in an arbitrary $\tilde \Psi$, moreover $\{\psi_{n,i}\}$ are linearly independent modulo ${\rm Im}\: d_{\Psi}$ 
because $\{[\psi_{n,i}]| i{=}7,{\dots},11\}$ is the basis in $\oplus_{{>}0}H^2_{(i)}(\mathfrak g, \mathfrak g)$.

The variety $V_n$ can be regarded as a subvariety of the affine variety ${\mathcal N}_n$ of $n$-dimensional 
nilpotent Lie algebras.
There exists a $GL_n$-action on ${\mathcal N}_n$ by basis changes. 
Let us consider the set $V_n'=\{y \in {\mathcal N}_n| \exists \Psi \in V_n, \exists g \in GL_n, y=g\star \Psi \}$.
\begin{remark}
In fact the Zariski closure of $V_n'$ in ${\mathcal N}_n$ coincides with one of irreducible components of 
${\mathcal N}_n$ that were discussed by Yu.Khakimdjanov in \cite{Khakim}.  
\end{remark}

The $GL_n$-action on $V_n'=\{y \in {\mathcal N}_n| \exists x \in V_n, \exists g \in GL_n, y=g\star x \}$
has singularities. Namely by the lemma \ref{isom_lemma} 
the stabilizer $(GL_n)_{\tilde \Psi}$ of a point $\tilde \Psi \in {\mathbb K}^5 \subset V_n$ 
coincides with the stabilizer $(G_n)_{\tilde \Psi}$ and we have 
$$\dim (GL_n)_{\tilde \Psi}= \left\{ \begin{array}{l} n{+}2, \; {\rm if} \;\: \tilde \Psi \ne 0;\\
n{+}3, \; {\rm if} \;\: \tilde \Psi = 0. \end{array} \right.$$
\begin{remark}
If $\tilde \Psi =0$ then one have to add the operator $\sum_{i{=}1}^{n}i e_i\otimes e^i$ to 
(\ref{basis_yadra_stab}) in order to get a basis of the Lie algebra of the stabilizer
$(GL_n)_{\tilde \Psi}$.
\end{remark}
Moreover, even for non-zero points $\tilde \Psi$ there are some particularities. Let 
$\tilde \Psi$ be a generic point in ${\mathbb K}^5$, i.e. $x_i \ne 0, \forall i$, then
$$
(GL_n)_{\tilde \Psi}=(G_n)_{\tilde \Psi}=(UT_n)_{\tilde \Psi}.
$$ 
But if ${\mathbb K}={\mathbb C}$ then for instance for $\tilde \psi=(1,0,0,0,0)$ we have
$$
(GL_n)_{\tilde \Psi}=(G_n)_{\tilde \Psi}={\mathbb Z}_{n{-}11}\ltimes (UT_n)_{\tilde \Psi},
$$ 
where ${\mathbb Z}_{n{-}11}$ stands for the group of the roots of the unit:
$${\mathbb Z}_{n{-}11}=\left\{ \alpha \in {\mathbb C}\;|\; \alpha^{n{-}11}=1 \right\}.$$

It is interesting here to remark some parallels between our discussions and the theory 
of non-singular deformations, that was introduced in \cite{FiFu}.  
\begin{definition}[\cite{FiFu}]
Let $\mathfrak{g}$ be a Lie algebra with the comutator $[,]$.
Consider a formal one-parameter deformation
$$[x,y]_t=[x,y]+\sum\limits_{k\ge 1}\alpha_k(x,y)t^k$$
of $\mathfrak{g}$. A  deformation is called non-singular if there exists a formal
one-parameter family of linear transformations
$$\varphi_t(x)=x+\sum\limits_{l\ge 1}\beta_k(x)t^l$$
of $\mathfrak{g}$ and a formal (not necessarily invertible) parameter
change $u=u(t)$ which transform the deformation $[x,y]_t$ into a deformation
$${[x,y]'}_u=[x,y]+\sum\limits_{k\ge 1}{\alpha '}_k(x,y)u^k,
\quad \quad \quad \varphi_t^{{-}1}[\varphi_t(x),\varphi_t(y)]_t={[x,y]'}_{u(t)}$$
with the cocycle ${\alpha '}_1 \in C^2(\mathfrak{g}, \mathfrak{g})$
being not cohomologous to $0$. Otherwise the deformation is called singular.
\end{definition}
We have already remarked that one can associate to an arbitrary ${\mathbb Z}_{{>}0}$-deformation
$\Psi=\Psi_1+\Psi_2+\Psi_3+\dots$ of ${\mathcal V}_n$ the following
one-parametric deformation:
$$[x,y]_t^{\Psi}=[x,y]+t\Psi_1(x,y)+t^2\Psi_2(x,y)+\dots+t^k\Psi_k(x,y)+\dots.$$
One can also remark that in our case a formal deformation
$[x,y]_t^{\Psi}$ is in fact a polynomial on $t$.
\begin{proposition}
Let fix the isomorphism $O(UT_n,V_n)={\mathbb K}^5, n\ge 16$, then 
a formal deformation $[x,y]_t^{\Psi}$ corresponding to some ${\mathbb Z}_{{>}0}$-deformation
$\Psi$ of ${\mathcal V}_n$ is non-singular if and only if its orbit $O_{\Psi}$ 
belongs to the union $\bigcup \limits_{i=1}^5  Ox_i$ of coordinate lines $Ox_i$ in ${\mathbb K}^5$.
\end{proposition} 
Acting by corresponding $\varphi_t=id+t\varphi_1+t^2\varphi_2+\dots$ we reduce our problem to a parameter change
in $$\Psi_t=t^{n{-}11} x_1\psi_{n,11}+t^{n{-}10} x_2\psi_{n,10}+\dots+t^{n{-}7} x_5\psi_{n,7}.$$ 
If $x_1 \ne 0$ then a parameter change $u(t)=t^{n{-}11}$ is possible iff $x_2=x_3=x_4=x_5=0$ the other cases
are treated in the same way.

\begin{remark}
If $n\ge 14$ one can associate to an arbitrary deformation of the form    
$\Psi=\sum_{i=1}^5 x_i\psi_{n,12{-}i}$ a linear one-parameter deformation $[,]_{\tau}'$
of ${\mathcal V}_n$:
$$[x,y]_{\tau}'=[x,y]+{\tau}\Psi(x,y).$$ 
It follows from $[\Psi, \Psi]=0$ in $C^3({\mathcal V}_n,{\mathcal V}_n)$.  
Evidently a linear one-parameter deformation $[,]_{\tau}'$ is non-singular 
iff $\Psi \ne 0$.
\end{remark}

\section{Symplectic structures and one-dimensional central extensions}

First of all we are going to calculate $H^2(\mathfrak g)$ of an arbitrary
${\mathbb Z}_{{>}0}$-deformation $\mathfrak g=({\mathcal V}_n,[,]+\Psi)$.
As we remarked in the section \ref{cohomology} the canonical basis $e_1,e_2,\dots,e_n$ of 
$\mathfrak g=({\mathcal V}_n,[,]+\Psi)$
defines the filtration $L$ of $\mathfrak g$: 
$$\mathfrak g=L^1 \mathfrak g \supset L^2 \mathfrak g {=}Span ( e_{2}, \dots, e_{n} )\supset \dots \supset
L^{n}\mathfrak g{=}Span (e_n) \supset \{0\},$$
and $L$ defines
the filtration $\tilde L$ of the cochain complex $(C^*(\mathfrak g),d)$.
For the simplicity one can assume that $\Psi=x_1\psi_{n,11}+\dots+x_5\psi_{n,7}$.
Let us consider the corresponding to $\tilde L$ spectral sequence $E_r^{p,q}$ that converges to $H^*(\mathfrak g)$.
We recall also that 
$$E^{p,q}_1 = H^{p{+}q}_{(-p)}({\rm gr}_L \mathfrak g)= H^{p{+}q}_{(-p)}({\mathcal V}_n).$$  
We recall that the space $H^2({\mathcal V}_n)$ is $3$-dimensional
and it is spanned
by classes 
$$g_5=[e^2 {\wedge} e^3],\;
g_7=[e^2 {\wedge} e^5 - 3 e^3 {\wedge} e^4], \;
[\Omega_{n{+}1}]{=}\frac{1}{2}[\sum\limits_
{i{+}j{=}n{+}1}(j{-}i)
e^i{\wedge}e^j]$$ of weights $5$, $7$, $n{+}1$ respectively.
It is evident that the classes $[g_5]$ and $[g_7]$ survive 
to the term $E_{\infty}$ because 
$d_{\tilde \Psi}e^i=d e^i, i=1,\dots,7,\; \forall \tilde \Psi$.
And the question is: does $[\Omega_{n{+}1}] \in E_1^{{-}2k{-}1,2k{+}3}$
survive to $E_{\infty}$ ?

\begin{proposition} Let $\mathfrak g=({\mathcal V}_n,[,]+\tilde \Psi)$ be a Lie algebra defined by
the following deformation of ${\mathcal V}_n$:  
$$\tilde \Psi=\tilde \Psi_{n{-}11}+\dots+\tilde \Psi_{n{-}7}=x_1\psi_{n,11}+\dots+x_5\psi_{n,7}.$$
Then in the spectral sequence $E_r^{p,q}$ the following properties
hold on:

1) $d_r \equiv 0,\: r=1,\dots,n{-}12$;

2) $E_r^{p,q}=E_1^{p,q}, \: r=2,\dots,n{-}11$; 

3) $d_{n{-}11}( \Omega_{n{+}1}) = 0$ if and only if $x_1=0$.

4) If $x_1=0$ then $d_r ( \Omega_{n{+}1}) = 0\;\; \forall r \ge 1$, i.e
the class $[\Omega_{n{+}1}]$ survive to the term $E_{\infty}$.
\end{proposition}

The differential $d_{\Psi}$ of the cochain complex $(C^*(\mathfrak{g}^*),d)$ has the form:
$$d=d_0+\Psi_{n{-}11}^*+\dots+\Psi_{n{-}7}^*=d_0+x_1\psi_{n,11}^*+\dots+x_5\psi_{n,7}^*,$$
where by $\Psi_{i}^*$ and $\psi_{n,j}^*$ we denote the dual applications to
$\Psi_{i}$ and $\psi_{n,j}$ respectively and $d_0$ stands for the differential 
of $(C^*({\mathcal V}_n),d_0)$. 
The items 1) and 2) of the proposition are evident and
$$d_{n{-}11}([\Omega_{n{+}1}])=[\Psi^*_{n{-}11}(\Omega_{n{+}1})]=x_1[\psi_{n,11}^*(\Omega_{n{+}1})],$$
On another hand one can compute $\psi_{n,11}^*(\Omega_{n{+}1})$ explicitely:
\begin{equation}
\begin{split}
-\psi_{n,11}^*(\Omega_{n{+}1})= 
-\psi_{n,11}^*\left((n{-}1)e^1\wedge e^n+(n{-}3)e^2\wedge e^{n{-}1}+\dots+(n{-}9)e^5\wedge e^{n{-}4}\right)=\\
(n{-}1)e^1\wedge 
(Z_1(n{-}4)e^2 {\wedge} e^9{+}Z_2(n{-}4)e^3 {\wedge} e^8
{+}Z_3(n{-}4)e^4 {\wedge} e^7{+}Z_4(n{-}4)e^5 {\wedge} e^6)+\\
(n{-}3)e^2\wedge (Q_1(n{-}4)e^2 {\wedge} e^8+Q_2(n{-}4)e^3 {\wedge} e^7
+Q_3(n{-}4)e^4 {\wedge} e^6)+\\
+(n{-}5)e^3 \wedge (P_1(n{-}4)e^2 {\wedge} e^7+P_2(n{-}4)e^3 {\wedge} e^6
+P_3(n{-}4)e^4 {\wedge} e^5)+\\
+ (n{-}7)e^4 \wedge \frac{n{-}5}{2} (e^2 {\wedge} e^6 - 2 e^3 {\wedge} e^5)+
+(n{-}9)e^5 \wedge (e^2 {\wedge} e^5 - 3 e^3 {\wedge} e^4)
\end{split}
\end{equation}
As it was shown in the Section \ref{Computations} this element is cohomologous
to 
$$-\frac{1}{5544}(n{-}12)(n^2{-}12n{+}59)(n^2{-}21n{+}116)
[2e^2{\wedge} e^3{\wedge} e^7{-}
5e^2{\wedge} e^4{\wedge} e^6{+}20e^3{\wedge} e^4{\wedge} e^5]$$
and it defines a non-trivial element in $H^3_{(12)}({\mathcal V}_n)$
for $n > 12$. 

We can finish the proof by remark  that
$$E_r^{{-}n{-}1+r,n{-}r{+}4}=E_1^{{-}n{-}1+r,n{-}r{+}4}=  H^3_{(n{+}1{-}r)}({\mathcal V}_n)=0, \; r>n-11,$$
because $H^3_{(i)}({\mathcal V}_n)=0$ for all  $i < 12$. Hence all other differentials of the spectral sequence
$$d_r:E_r^{{-}n{-}1,n{+}3} \to  E_r^{{-}n{-}1+r,n{-}r{+}4},\; r>n{-}11$$ are trivial.

\begin{corollary}
Let $\mathfrak g=({\mathcal V}_{2k},[,]+ \Psi)$ be a  ${\mathbb Z}_{{>}0}$-deformation of ${\mathcal V}_{2k}$
and $(x_1,x_2,\dots,x_5)$ be the set of homogeneous coordinates of the $UT_n$-orbit $O_{\Psi}$.
$\mathfrak g$ admits a symplectic structure if and only if $x_1=0$.
\end{corollary}  

\begin{corollary}
If the orbit $O_{\Psi}$ of a ${\mathbb Z}_{{>}0}$-deformation $\mathfrak g=({\mathcal V}_n,[,]+ \Psi$ has a non-zero
first coordinate $x_1 \ne 0$ then:
$$\dim H^1(\mathfrak g)=\dim H^2(\mathfrak g)=2.$$
\end{corollary}

Hence for a generic deformations $\mathfrak g$ of ${\mathcal V}_n$ the Dixmier inequalities (see \cite{D})
$\dim H^i(\mathfrak g) \ge 2$ are sharp for $i=1,2$.

Recall that a one-dimensional central extension of a Lie algebra $\mathfrak{g}$ is an exact sequence
\begin{equation}
\label{exactseq}
0 \to \mathbb K \to \tilde {\mathfrak g}  \to \mathfrak{g} \to 0 
\end{equation}
of Lie algebras and their homomorphisms, in which the image of the homomorphism
${\mathbb K} \to \tilde {\mathfrak{g}}$
is contained in the center of the Lie algebra $\mathfrak{g}$.
To the cocycle $c \in \Lambda^2(\mathfrak{g}^*)$ corresponds the extension
$$
0 \to {\mathbb K} \to {\mathbb K} \oplus \mathfrak{g} \to \mathfrak{g} \to 0
$$
where the Lie bracket in $\mathbb{K} \oplus \mathfrak{g}$ is defined by the formula
$$
[(\lambda,g), (\mu, h)]=(c(g,h), [g,h]).
$$
It can be checked directly that the Jacobi identity for this Lie bracket is equivalent to $c$
being cocycle and that to cohomologous cocycles correspond equivalent (in a obvious sense)
extensions.

Now let $\tilde{\mathfrak{g}}$ be a filiform Lie algebra, it has one-dimensional center $Z(\tilde{\mathfrak{g}})$
and we have the following one-dimensional central extension:
$$0 \to \mathbb K=Z(\tilde{\mathfrak{g}}) \to \tilde {\mathfrak g}  \to \tilde {\mathfrak g}/ Z(\tilde{\mathfrak{g}})  \to 0 $$
where the quotient Lie algebra  $\tilde {\mathfrak g}/ Z(\tilde{\mathfrak{g}})$ is also filiform Lie algebra.
Moreover, let $e_1,\dots,e_{n{+}1}$ be some adapted basis of $\tilde{\mathfrak g}$ ($Z(\tilde{\mathfrak g})=Span(e_{n{+}1})$) 
and $e^1,\dots, e^{n{+}1} \in {\mathfrak g}^*$ its dual basis.
Then the forms $e^1,\dots, e^n$ can be regarded as the dual basis to the basis 
$e_1+Z(\tilde{\mathfrak g}),\dots, e_n+Z(\tilde{\mathfrak g})$ of $\tilde {\mathfrak g}/ Z(\tilde{\mathfrak{g}})$
and the cocycle $\Omega_{n{+}1}=de^{n{+}1}$ determines this one-dimensional central extension.

\begin{proposition}
Let $\mathfrak g$ be a filiform Lie algebra 
 and let its center $Z({\mathfrak{g}})$ be spanned by some $\xi \in Z({\mathfrak{g}})$.
Then $\tilde {\mathfrak g}$ taken from a one-dimensional central extension 
$$0 \to \mathbb K \to \tilde {\mathfrak g}  \to \mathfrak{g} \to 0 $$
with cocycle $c \in \Lambda^2(\mathfrak g)$ 
 is a filiform Lie algebra if and only if
the restricted function $f( \cdot)=c( \cdot, \xi)$ is non-trivial in $\mathfrak{g}^*$.
\end{proposition}

\begin{corollary}
Let $\mathfrak g$ be a symplectic filiform Lie algebra and $\omega$ its symplectic form,
then the one-dimensional central extension $\tilde{\mathfrak g}$ with the cocycle $c=\omega$ will be
also a filiform Lie algebra.
\end{corollary}

\begin{theorem}
Let $\tilde{\mathfrak g} =({\mathcal V}_{2k{+}1},[,]+\tilde \Psi)$ be  
${\mathbb Z}_{{>}0}$-filtered deformation of ${\mathcal V}_{2k{+}1}$, where
$$\tilde \Psi=x_1\psi_{2k{+}1,11}+\dots+x_5\psi_{2k{+}1,7}.$$ 
Then $\tilde{\mathfrak g}$ can be represented as a one-dimensional central extension
of ${\mathbb Z}_{{>}0}$-filtered deformation ${\mathfrak g}_X =({\mathcal V}_{2k},[,]+\Phi)$ of ${\mathcal V}_{2k}$, where
$$\Phi=x_1\psi_{2k,10}+x_2\psi_{2k,9}+x_3\psi_{2k,8}+x_4\psi_{2k,7}, \; X=(x_1,x-2,x_3,x_4).$$ 
The cocycle $\Omega_{X,x_5}$ that determines
this one-dimensional central extension is equal to
\begin{equation}
\begin{split}
\label{can_sympl}
\Omega_{X,x_5}=\frac{1}{2}\sum\limits_{i{+}j{=}2k{+}1}(j{-}i)e^i{\wedge}e^j{+}
x_1\Omega_{2k{,}11}{+}x_2\Omega_{2k{,}10}{+}x_3\Omega_{2k{,}9}{+}x_4\Omega_{2k{,}8}{+}x_5\Omega_{2k{,}7},\\
\Omega_{2k{,}11}=
Z_1(2k{-}3)e^2{\wedge}e^9{+}Z_2(2k{-}3)e^3{\wedge}e^8{+}Z_3(2k{-}3)e^4{\wedge}e^7{+}Z_4(2k{-}3)e^5{\wedge}e^6,\\
\Omega_{2k{,}10}=Q_1(2k{-}2)e^2{\wedge}e^8+Q_2(2k{-}2)e^3{\wedge}e^7+Q_3(2k{-}2)e^4{\wedge}e^6,\\
\Omega_{2k{,}9}=P_1(2k{-}1)e^2{\wedge}e^7+P_2(2k{-}1)e^3{\wedge}e^6+P_3(2k{-}1)e^4{\wedge}e^5,\\
\Omega_{2k{,}8}=\frac{2k{-}1}{2}(e^2{\wedge}e^6-2e^3{\wedge}e^5),\;\; \Omega_{2k{,}7}=e^2{\wedge}e^5-3e^3{\wedge}e^4.
\end{split}
\end{equation}
\end{theorem}

\begin{definition}
Let $\mathfrak g, \tilde {\mathfrak g}$ be two symplectic Lie algebras,
$\omega_{\mathfrak g}, \omega_{\tilde {\mathfrak g}}$ be corresponding symplectic structures.
A Lie algebras isomorphism $f: \mathfrak g \to  \tilde {\mathfrak g}$ is called
a symplecto-isomorphism $f: (\mathfrak g,\omega_{\mathfrak g}) \to  (\tilde {\mathfrak g}, \omega_{\tilde {\mathfrak g}})$
if and only if 
$\omega_{\mathfrak g}=f^*(\omega_{\tilde {\mathfrak g}})$.
\end{definition}

\begin{theorem}
1) Let $\mathfrak g$ be a symplectic filtered deformation of ${\mathcal V}_{2k}$,
$\omega_{\mathfrak g}$ its symplectic structure. 
Then there exists a vector $X=(x_1,x_2,x_3,x_4) \in {\mathbb K}^4$,  $x_5 \in {\mathbb K}$ 
such that the pair  $(\mathfrak g,\omega_{\mathfrak g})$ is symplecto-isomorphic
to $({\mathfrak g}_X, \Omega_{X{,}x_5})$, where  
${\mathfrak g}_X$ is one of the Lie algebras
from the table below 
$$
\begin{tabular}{|c|c|}
\hline
&\\[-10pt]
$i+j$ & ${\mathfrak g}_X$: commutating relations $[e_i,e_j], \; i<j$ \\
&\\[-10pt]
\hline
&\\[-10pt]
$3 \le i+j \le 6$ & 
$[e_i,e_j]=(j-i)e_{i+j}$ \\
&\\[-10pt]
\hline
&\\[-10pt]
$7$ & \begin{tabular}{c}
$[e_1,e_6]=5e_7$, \\
$[e_2,e_5]=3e_7+x_1e_{2k{-}3}+x_2e_{2k{-}2}+x_3e_{2k{-}1}
+x_4e_{n}$,\\
$[e_3,e_4]=e_7-3x_1e_{2k{-}3}-3x_2e_{2k{-}2}-3x_3e_{2k{-}1}
-3x_4e_{n}$,\\
\end{tabular}\\
&\\[-10pt]
\hline
&\\[-10pt]
$8$ & \begin{tabular}{c}
$[e_1,e_7]=6e_8$, \\
$[e_2,e_6]=4e_8+x_1 \frac{(2k{-}4)}{2}e_{2k{-}2}+
x_2\frac{(2k{-}3)}{2}e_{2k{-}1}+x_3 \frac{(2k{-}2)}{2}e_{2k}$,\\
$[e_3,e_5]=2e_8+x_1(2k{-}4)e_{2k{-}2}+x_2(2k{-}3)e_{2k{-}1}+
x_3(2k{-}2)e_{2k}$;\\
\end{tabular}\\
&\\[-10pt]
\hline
&\\[-10pt]
$9$ & \begin{tabular}{c}
$[e_1,e_8]=7e_9$, \\
$[e_2,e_7]=5e_9+x_1
\frac{(5\binom{2k{-}3}{2}{+}3\binom{2k{-}3}{1}{-}6)}{21}e_{2k{-}1}+
x_2
\frac{(5\binom{2k{-}2}{2}{+}3\binom{2k{-}2}{1}{-}6)}{21}e_{2k}$,\\
$[e_3,e_6]=3e_9-x_1
\frac{(4\binom{2k{-}3}{2}{+}15\binom{2k{-}3}{1}{-}30)}{21}e_{2k{-}1}-
x_2
\frac{(4\binom{2k{-}2}{2}{+}15\binom{2k{-}2}{1}{-}30)}{21}e_{2k}$,\\
$[e_4,e_5]=e_9-x_1
\frac{(13\binom{2k{-}3}{2}{-}30\binom{2k{-}3}{1}{+}60)}{21}e_{2k{-}1}-
x_2
\frac{(13\binom{2k{-}2}{2}{-}30\binom{2k{-}2}{1}{+}60)}{21}e_{2k{-}1}$,\\
\end{tabular}\\
&\\[-10pt]
\hline
&\\[-10pt]
$10$ & \begin{tabular}{c}
$[e_1,e_9]=8e_{10}$, \\
$[e_2,e_8]=6e_{10}+x_1
\frac{(3\binom{2k{-}2}{3}{+}4\binom{2k{-}2}{2}{-}4\binom{2k{-}2}{1})}{28}
e_{2k}$,\\
$[e_3,e_7]=4e_{10}+x_1
\frac{(\binom{2k{-}2}{3}{-}8\binom{2k{-}2}{2}{+}8\binom{n{-}2}{1})}{14}
e_{2k}$,\\
$[e_4,e_6]=2e_{10}+x_1
\frac{({-}13\binom{2k{-}2}{3}{+}20\binom{2k{-}2}{2}{-}20\binom{2k{-}2}{1})}{28}
e_{n}$,\\
\end{tabular}\\
&\\[-10pt]
\hline
&\\[-10pt]
$11 \le i+j \le 2k$ & 
$[e_i,e_j]=(j-i)e_{i+j}$ \\
\hline
\end{tabular}$$
and its corresponding symplectic form $\Omega_{X,x_5}$ is equal to 
$$\Omega_{X,x_5}=\frac{1}{2}\sum\limits_{i{+}j{=}2k{+}1}(j{-}i)e^i{\wedge}e^j{+}
x_1\Omega_{2k{,}11}{+}x_2\Omega_{2k{,}10}{+}x_3\Omega_{2k{,}9}{+}x_4\Omega_{2k{,}8}{+}x_5\Omega_{2k{,}7}.$$
2) A pair $(\mathfrak g_X, \Omega_{X,x_5})$ is symplecto-isomorphic
to $(\mathfrak g_Y, \Omega_{Y,y_5})$
if and only if
there exist an $\alpha \in {\mathbb K}^*$ such that
$$
y_1=\alpha^{n{-}11} x_1,y_2=\alpha^{n{-}10} x_2,
y_3=\alpha^{n{-}9} x_3,y_4=\alpha^{n{-}8} x_4, y_5=\alpha^{n{-}7} x_5. 
$$ 
\end{theorem}
\begin{remark}
The previous theorem shows that all ${\mathbb Z}_{{>}0}$-deformations $\tilde {\mathfrak g}$ of ${\mathcal V}_{2k{+}1}$
are contact Lie algebras and gives their complete classification. 
An arbitrary symplectic ${\mathbb Z}_{{>}0}$-deformations $\mathfrak g$ of ${\mathcal V}_{2k}$ can be
obtained as quotient Lie algebra $\tilde {\mathfrak g}/Z(\tilde {\mathfrak g})$. This method
of classification of symplecto-isomorphism classes of low-dimensional filiform
Lie algebras was considered in \cite{GJKh2}.
\end{remark}

Taking rational coordinates $(x_1,x_2,\dots,x_5)$ in the table above
one will get a nilpotent Lie algebra ${\mathfrak g}_X$ with rational structure constants 
and hence due to the Malcev theorem (\cite{Mal}) the corresponding simply connected nilpotent
Lie group $G$ has a cocompact lattice $\Gamma$. Thus one will get
a family of examples of symplectic nilmanifolds $M=G/\Gamma$.

\end{document}